\documentclass{amsart}
\usepackage{graphicx}
\usepackage{amssymb}
\usepackage{epstopdf}
\usepackage{nicefrac}
\usepackage{enumerate}
\usepackage{hyperref}
\usepackage{float}
\usepackage{color}
\usepackage{pst-all}
\usepackage{tabularx}

\input xy
\xyoption {all}

\newtheorem{thm}{THEOREM}[section]

\newtheorem{cor}[thm]{COROLLARY}
\newtheorem{defn}[thm]{DEFINITION}

\newtheorem{lemma}[thm]{LEMMA}
\newtheorem{prob}[thm]{PROBLEM}
\newtheorem{prop}[thm]{PROPOSITION}
\newtheorem{quest}[thm]{QUESTION}
\newtheorem{remark}[thm]{REMARK}

\newcommand{\ds}{\displaystyle}

\newcommand{\cA}{{\mathcal A}}
\newcommand{\cB}{{\mathcal B}}

\newcommand{\cG}{{\mathcal G}}
\newcommand{\cH}{{\mathcal H}}
\newcommand{\cK}{{\mathcal K}}
\newcommand{\cN}{{\mathcal N}}
\newcommand{\cP}{{\mathcal P}}
\newcommand{\cQ}{{\mathcal Q}}
\newcommand{\cS}{{\mathcal S}}
\newcommand{\cY}{{\mathcal Y}}
\newcommand{\dF}{d_{\F}} 
\newcommand{\e}{{\epsilon}}
\newcommand{\eF}{{\epsilon_{\F}}}    
\newcommand{\F}{{\mathcal F}}
\newcommand{\FfM}{{\mathcal F}_{\fM}} 
\newcommand{\fG}{\mathfrak{G}}
\newcommand{\fM}{{\mathfrak{M}}}
\newcommand{\FP}{{\mathcal F}_{\cP}} 
 
\newcommand{\fT}{{\mathfrak{T}}}
\newcommand{\fX}{{\mathfrak{X}}}
\newcommand{\lF}{{\lambda_{\mathcal F}}}
\newcommand{\mH}{{\mathbb H}}
\newcommand{\mN}{{\mathbb N}}
\newcommand{\mQ}{{\mathbb Q}}
\newcommand{\mR}{{\mathbb R}}
\newcommand{\mS}{{\mathbb S}}
\newcommand{\mT}{{\mathbb T}}
\newcommand{\mZ}{{\mathbb Z}}
\newcommand{\wtzM}{\widetilde{M_0}}

\begin{document}

\title{Pro-groups and generalizations of a theorem of Bing}

\thanks{2010 {\it Mathematics Subject Classification}. Primary 52C23, 57R05, 54F15, 37B45; Secondary 53C12, 57N55 }

\author{Alex Clark}
\address{Alex Clark, School of Mathematical Sciences, Queen Mary University of London, London E1 4NS, UK}
\email{alex.clark@qmul.ac.uk}

\author{Steven Hurder}
\address{Steven Hurder, Department of Mathematics, University of Illinois at Chicago, 322 SEO (m/c 249), 851 S. Morgan Street, Chicago, IL 60607-7045}
\email{hurder@uic.edu}

\author{Olga Lukina}
\address{Olga Lukina, Department of Mathematics, University of Illinois at Chicago, 322 SEO (m/c 249), 851 S. Morgan Street, Chicago, IL 60607-7045}
\email{ollukina940@gmail.com}

 \thanks{Version date: November 1, 2018 }

\date{}



\begin{abstract}
 A matchbox manifold is a generalized lamination, which is a continuum whose path components define the leaves of a foliation of the space.
A matchbox manifold is $M$-like if it has the shape of a fixed topological space $M$. When $M$ is a closed manifold, in a previous work, the authors have shown that if   $\fM$ is a matchbox manifold which is $M$-like, then it is homeomorphic to a weak solenoid. In this work, we associate to a weak solenoid  a pro-group, whose pro-isomorphism class  is an invariant of the homeomorphism class of $\fM$.   We then show that   an $M$-like matchbox manifold is   homeomorphic to a weak solenoid  whose base manifold has fundamental group which is non co-Hopfian; that is, it admits a non-trivial   self-embedding of finite index. We include a collection of examples illustrating this conclusion.
\end{abstract}

\maketitle



\section{Introduction} \label{sec-intro}

A \emph{continuum} is a compact \emph{connected} metric space. From the very beginnings of the subject of Topology in the first half of the twentieth century, examples of continua with unusual properties were   a source of motivation for the development of the subject. For example,  the Warsaw circle, the Sierpinsky carpet \cite{Sierpinski1916}, and the 1-dimensional solenoids introduced by Vietoris \cite{Vietoris1927} and van Danzig \cite{vanDantzig1930} are all well-known examples of ``non-trivial'' continua. Continua frequently arise as invariant sets for dynamical systems, and moreover as quoted from Barge and Kennedy \cite{BargeKennedy1990}:
\begin{quote}
One can argue quite convincingly that continuum theory first arose from problems in dynamics even before there was a definition of a topological space.
\end{quote}

The study of the   invariant continua for a   dynamical system is a basic technique for the analysis of its global dynamics.  
There are many established approaches   for the study of continua, as   discussed in the book by Sam Nadler \cite{Nadler1992}. One of these is the following fundamental result of Hans Freudenthal:

\begin{thm}[Freudenthal \cite{Freudenthal1936}] \label{thm-Freudenthal}
Let    $\fM$  be a continuum, then  there exists a sequence of maps
\begin{equation}\label{eq-presentationFreud}
\cP = \{\, f_{\ell+1} \colon  X_{\ell+1} \to X_\ell \mid  \ell \geq 0\}
\end{equation}
where for each $\ell \geq 0$,  $X_{\ell}$ is a finite connected polyhedron and  each   map  $f_{\ell +1}$ is a continuous surjection, such that
$\fM$ is homeomorphic to an inverse limit, 
\begin{equation}\label{eq-presentation-top}
\Phi_{\cP} \colon \fM ~ \cong  ~  \cS_{\cP} \equiv  \varprojlim \,\{\, f_{\ell+1} \colon  X_{\ell+1} \to X_\ell \mid  \ell \geq 0\} ~ ,
\end{equation}
 The maps  $f_{\ell +1}$ in \eqref{eq-presentationFreud} are called the \emph{bonding maps} of the presentation $\cP$. 
\end{thm}
  The inverse limit is considered as a subset of the Cartesian product of the spaces $\{X_{\ell} \mid \ell \geq 0\}$ with the Tychonoff topology, and a point   is represented by a sequence $x = (x_{\ell})_{\ell \geq 0}$ where $x_{\ell} \in X_{\ell}$ and $f_{\ell +1}(x_{\ell+1}) = x_{\ell}$ for $\ell \geq 0$. 
A presentation   $\cP$ is said to be \emph{trivial} if  all of the maps $f_{\ell +1}$ are homeomorphisms, and thus  $\fM$ is homeomorphic to a finite connected polyhedron.

 Observe that given a presentation $\cP$ for a continuum $\fM$ as in \eqref{eq-presentation-top},     for all $\ell \geq 0$, there is a continuous surjection $\Pi_{\ell} \colon \fM \to X_{\ell}$ given by $\Pi_{\ell}((x_{\ell})_{\ell \geq 0}) = x_{\ell}$. It then follows that for any $\e > 0$, there exists some $\ell_{\e}$ such that  $\Pi_{\ell_\e}$ is an $\e$-mapping; that is, for all $y \in X_{\ell_\e}$ the inverse image set  $\Pi_{\ell_\e}^{-1}(y) \subset \fM$ has diameter less than $\e$. 
 The notion of  an $\e$-map  was introduced in 1928 by     Alexandroff     \cite{Alexandroff1928}: 
\begin{defn}    \label{def-Ylike}
 Let  $\fM$ be a continuum, $Y$ a topological space and $\e > 0$ a constant. A   map $f \colon \fM \to Y$ is  said to be an $\e$-map if $f$ is a continuous surjection, and for each point $y\in Y$,   the inverse image $f^{-1}(y)$   has diameter less than $\e$.   A continuum $\fM$ is said to be \emph{$Y$--like}, for some topological space $Y$, if for every   $\e>0$, there is an $\e$-map $f_\e \colon \fM \to Y$. 
 More generally, let $\cY$ be a collection of compact   topological spaces, then $\fM$ is $\cY$-like if for all $\e > 0$, there exists $Y_\e \in \cY$ and an $\e$-map $f_\e \colon \fM \to Y_\e$. 
 \end{defn}
 
 In this work, we consider the following special case of the above definition:
\begin{defn}    \label{def-Mlike}
 A continuum  $\fM$ is \emph{manifold-like}, if given any $\e > 0$,  there exists a closed manifold $M_{\e}$ and an $\e$-map
 $f \colon \fM \to M_{\e}$. For a closed manifold $M$, a continuum  $\fM$ is \emph{$M$-like}, if given any $\e > 0$,  there exists  an $\e$-map
 $f \colon \fM \to M$.
  \end{defn} 
 
 The following is then a natural problem to investigate:
  \begin{prob}\label{prob-Mlike}
For a closed $n$-manifold $M$, characterize  the  continua which  are $M$-like. 
 \end{prob}

 Problem~\ref{prob-Mlike} is too broad  to obtain definitive results, as there exist too many pathological constructions of $M$-like spaces.  The simplest case, where $M = [0,1]$ is an interval, has been extensively studied for many special classes of bonding maps, such as the    ``tent maps'', where many results are known for the inverse limit spaces. However, without some sort of   ``regularity'' restriction, even the continua modeled on bonding maps to $[0,1]$ are too wild to hope to classify. Thus,  one imposes additional assumptions on the continua considered, such as the following.
 Recall that  a topological  space $\fM$ is \emph{homogeneous} if for
every $x, y \in \fM$, there exists a homeomorphism $h \colon \fM \to \fM$ such that $h(x) = y$. 
A continuum $X$ is \emph{circle-like} if it is $M$-like, where $M = \mS^1$ is the circle. 
The following is a well-known result of R.H. Bing, which has inspired many generalizations:
\begin{thm} [Bing \cite{Bing1960}] \label{thm-bing}
Let $\fM$ be    a homogeneous, circle-like
continuum that contains an arc. Then either $\fM$ is homeomorphic to   $\mS^1$, or to    an inverse limit of proper finite coverings of   $\mS^1$.
\end{thm}

The    hypothesis in Theorem~\ref{thm-bing}  that $\fM$ is homogeneous and contains an arc implies that for every point $x \in \fM$,  there is an arc containing $x$. One generalization of this hypothesis    is suggested by the properties of continua arising    from the investigation  of the dynamical properties of smooth foliations of compact manifolds, for example as discussed in \cite{Hurder2014}. In this case, the invariant continua have the structure of foliated spaces,   in the sense of \cite{MS2006}. That is, for each $x \in \fM$, there is a neighborhood $x \in U_x \subset \fM$ where $U_x$ is homeomorphic to a product of an open subset of $\mR^n$, for some $n \geq 1$, with a closed subset $\fT_x$ of some Polish space.  If the transversal spaces $\fT_x$ are \emph{totally disconnected} for all $x \in \fM$, and are not singleton spaces, then we say that  $\fM$ is a \emph{matchbox manifold}, and $\fM$ then admits  a decomposition $\FfM$ into path-connected components  of constant dimension $n$ which are  the leaves of $\FfM$.   A matchbox manifold with $2$-dimensional leaves is a lamination by surfaces in the sense of Ghys  \cite{Ghys1999} and Lyubich and Minsky \cite{LM1997}, while Sullivan called them ``solenoidal spaces'' in \cite{Sullivan2014,Verjovsky2014}. The terminology ``matchbox manifold'' follows  the usage introduced in  \cite{AM1988,AO1991,AO1995}, and as used in the authors' works 
\cite{ClarkHurder2011,ClarkHurder2013,CHL2014,CHL2018a,CHL2018b,DHL2016a,DHL2016b,DHL2016c,HL2017,HL2018}
  which study this class of continua. 
  
In the authors' previous work \cite{CHL2018b}  we have shown:  
 \begin{thm}\cite[Corollary~1.6]{CHL2018b} \label{thm-weak}
Let $\fM$ be a    manifold-like matchbox manifold $\fM$. Then there exists a   presentation  $\cP = \{ p_{\ell +1} \colon M_{\ell +1} \to M_{\ell}  \mid  \ell \geq 0 \}$,  
 where each   bonding map  $p_{\ell +1}$ is a proper covering map, and a homeomorphism
 \begin{equation}\label{eq-ws}
\Phi_{\cP} \colon \fM ~  \cong  ~ \cS_{\cP} \equiv \varprojlim\,   \{ p_{\ell +1} \colon M_{\ell +1} \to M_{\ell}  \mid  \ell \geq 0 \} ~.
\end{equation}
\end{thm}
A continuum given by   the inverse limit of a sequence of proper covering maps is called a \emph{weak solenoid}, following the usage in the the works of  Rogers and Tollefson \cite{RogersJT1970,RT1971a,RT1971b,RT1972a}. The properties of these spaces are discussed further in Section~\ref{sec-solenoids}.

Note that a homeomorphism $\Phi_{\cP} \colon \fM \to \cS_{\cP}$ preserves path components, hence maps the leaves of the foliation $\FfM$ homeomorphically onto the leaves of the foliation $\FP$ of $\cS_{\cP}$.

    Mardesic and Segal discuss pointed shape expansions in \cite{MardesicSegal1982}, and show there is  well-defined notion of pro-homotopy groups {\it pro}-$\pi_k(\fM, z)$ of a space $\fM$ using pointed shape exapnsions, though it may depend on the choice of basepoint $z \in \fM$. 
For example, see Example~4,  \cite[Chapter II, Section 3.3]{MardesicSegal1982}. One of the main technical points of this work, is that the pro-homotopy groups are well-defined for a weak solenoid,   independent of the choice of basepoint $z \in \fM$. 

Here is our first result, which is a general property of weak solenoids, and can be considered a formal version of the ideas of Fokkink and Oversteegen in \cite{FO2002}.

 \begin{thm}\label{thm-MMprogroups}
Let $\fM$ be an $M$-like matchbox manifold, then the pro-group {\it pro}-$\pi_1(\fM, z)$  is independent of the choice of basepoint $z \in \fM$ up to isomorphism, and depends only on  the homeomorphism class of $\fM$. 
 \end{thm}
Pro-groups are discussed in   Section~\ref{sec-progroups}, then    the definition of {\it pro}-$\pi_1(\fM, z)$ is given    in Section~\ref{sec-progroupsWS} and we show that it is well-defined.
Then in Section~\ref{sec-progroupsMlike}, we calculate the pro-group associated to an $M$-like matchbox manifold, which leads to the following   generalization of Theorem~\ref{thm-bing} by Bing:
 \begin{thm}\label{thm-main1}
 Let $\fM$ be an $n$-dimensional  matchbox manifold which  is $M$-like, for a closed manifold $M$. 
 Then $\fM$   admits a presentation 
 $\ds \cP = \{\, q_{\ell+1} \colon  N_{\ell+1} \to N_{\ell} \mid  \ell \geq 0\}$, 
 where  each $N_{\ell}$ is a closed $n$-manifold whose  fundamental group is isomorphic to that of  $M$, and each $q_{\ell+1} \colon  N_{\ell +1} \to N_{\ell}$ is a proper covering map for $\ell \geq 0$, such that there is a homeomorphism  
  \begin{equation}\label{eq-scale1}
\Phi_{\cP} \colon \fM     ~  \cong ~ \cS_{\cP} \equiv     \varprojlim \,\{\, q_{\ell+1} \colon  N_{\ell +1} \to N_{\ell} \mid  \ell \geq 0\} ~ .  
\end{equation}
 \end{thm}

   The conclusion of  Theorem~\ref{thm-main1} can be strengthened if we impose additional hypotheses on the target space $M$.
Recall that a finite $CW$-complex  $Y$ is \emph{aspherical} if it is connected  and its universal covering space  is  contractible.     The \emph{Borel Conjecture} is that if  $Y_1$ and $Y_2$ are homotopy equivalent, aspherical closed manifolds, then a homotopy equivalence between $Y_1$ and $Y_2$ is homotopic to a  homeomorphism between $Y_1$ and $Y_2$.    The Borel Conjecture has been proven for many classes of aspherical manifolds:
\begin{enumerate}
\item the torus $\mT^n$ for all $n \geq 1$,
\item all  infra-nilmanifolds,
\item   closed Riemannian manifolds $Y$ with negative sectional curvatures,
\item   closed Riemannian manifolds $Y$ of dimension $n \ne 3,4$ with non-positive sectional curvatures. 
\end{enumerate}
 The above  list is not exhaustive, and only cases (1) and (2) are used in this paper. The history and   status of the Borel Conjecture is discussed, for example, in the surveys of Davis \cite{Davis2012} and L\"{u}ck  \cite{Luck2012}. 
 
  A manifold for which all finite coverings satisfy the Borel Conjecture was called \emph{strongly Borel} in the authors work \cite{CHL2018a}. Note that the class of ``D-like'' spaces introduced in \cite[Section~5]{RT1972b} consists of manifolds which are strongly Borel.
  Here is our second main result, which is a further generalization of Theorem~\ref{thm-bing} by Bing:

 \begin{thm}\label{thm-main2}
 Let $M$ be a closed aspherical $n$-manifold $M$ which satisfies the Borel Conjecture. If $\fM$ is a matchbox manifold which  is $M$-like, 
  then $\fM$   admits a presentation 
 $\ds \cP = \{\, q_{\ell+1} \colon  M \to M \mid  \ell \geq 0\}$, 
where    each $q_{\ell+1} \colon  M \to M$ is a proper covering map for $\ell \geq 0$,   such that there is a homeomorphism  
  \begin{equation}\label{eq-scale2}
\Phi_{\cP}   \colon \fM  ~  \cong   ~ \varprojlim \,\{\, q_{\ell+1} \colon  M \to M \mid  \ell \geq 0\} ~ .  
\end{equation}
 \end{thm}

 As an example, note that the $n$-torus $\mT^n$ satisfies the hypotheses of Theorem~\ref{thm-main2}, so we   conclude:
 \begin{cor}\label{cor-toral}
Let $\fM$ be a matchbox manifold which  is $\mT^n$-like, then it admits a presentation as in \eqref{eq-scale2} for $M = \mT^n$.
 \end{cor}

Recall that a group $H$ is said to be \emph{(finitely) non co-Hopfian}, or just non co-Hopfian in this work,  if there is an injective homomorphism $\iota \colon H \to H$ whose image has finite index, but is not onto. That is, $H$ contains a subgroup of finite index isomorphic to itself.  The class of non co-Hopfian groups has been extensively studied, as discussed further in Section~\ref{sec-examples},   providing many examples. The general philosophy is that the group $H$ should be ``essentially'' nilpotent.  A complete solution of Problem~\ref{prob-Mlike} in the case of matchbox manifolds then requires a classification of the finitely-generated, non co-Hopfian groups. The works of van~Limbeek \cite{vanLimbeek2018,vanLimbeek2017} give some indications of the available techniques for the study of this problem, as well as some of the difficulties   it presents.

The proofs of Theorems~\ref{thm-main1} and \ref{thm-main2} are given in Section~\ref{sec-proofs}.   Finally, Section~\ref{sec-examples}     gives   collections of examples of weak solenoids  and non co-Hopfian groups  illustrating the conclusions of  theorems, and showing that stronger results are not possible without further assumptions.

\section{Weak solenoids}\label{sec-solenoids}

  A  \emph{weak solenoid} $\cS_\cP$ of dimension $n \geq 1$, is the inverse limit space of a sequence 
  \begin{eqnarray} 
\cS_{\cP} & \equiv &  \varprojlim  ~ \{ p_{\ell +1} \colon M_{\ell +1} \to M_{\ell}\} \label{eq-presentationinvlim} \\
& = &  \{(x_0, x_1, \ldots ) \in \cS_{\cP}  \mid p_{\ell +1 }(x_{\ell + 1}) =  x_{\ell} ~ {\rm for ~ all} ~ \ell \geq 0 ~\} ~ \subset \prod_{\ell \geq 0} ~ M_{\ell} ~   \nonumber
\end{eqnarray}
   where each $M_\ell$ is a compact connected $n$-dimensional manifold without boundary, and  each  $p_{\ell +1 }$ is a   covering map  of finite degree greater than one.
The set $\cS_{\cP}$ is given  the relative  topology, induced from the product (Tychonoff) topology, so that $\cS_{\cP}$ is itself compact and connected.   The collection  
\begin{equation}\label{eq-presentation}
\cP = \{p_{\ell +1} \colon M_{\ell +1} \to M_{\ell} \mid \ell \geq 0\},
\end{equation}
 is called a \emph{presentation} for $\cS_\cP$.
 For example, a \emph{Vietoris solenoid} \cite{vanDantzig1930,Vietoris1927}  is   a $1$-dimensional    solenoid, where each $M_{\ell}$ is a circle, as arises in the conclusion of Theorem~\ref{thm-bing}.  

 Let $\cS_\cP$ be a weak solenoid with presentation $\cP$ as in \eqref{eq-presentation}. For each $\ell \geq 1$, the composition
\begin{equation}\label{eq-coverings}
p^0_{\ell} = p_{1} \circ \cdots \circ p_{\ell -1} \circ p_{\ell} \colon M_{\ell} \to M_0 ~ 
\end{equation}
is a finite-to-one covering map of the base manifold $M_0$.
  For each $\ell \geq 0$, projection onto the $\ell$-th factor in the product $\ds \prod_{\ell \geq 0} ~ M_{\ell}$ in \eqref{eq-presentationinvlim} yields a    surjective map denoted by $\Pi_{\ell}^{\cP} \colon \cS_{\cP}  \to M_{\ell}$. Then for  $\ell \geq 1$, we have  
$ \Pi_0^{\cP}  =   p^0_{\ell} \circ \Pi_{\ell}^{\cP}   \colon \cS_{\cP} \to M_0$. 
  
McCord showed in \cite{McCord1965} that for each $\ell > 0$, the map $\Pi_{\ell}^{\cP}  \colon \cS_{\cP} \to M_{\ell}$  is a fibration over $M_{\ell}$.  
For a choice of a basepoint $x_0 \in M_0$, the assumption  that the fibers of each map $p_{\ell+1}$ have cardinality at least $2$ implies the  fiber  $\fX_0 = (\Pi_0^{\cP})^{-1}(x_0)$ is a Cantor set. 
McCord also observed that   $\cS_\cP$ has a local product structure, which  implies that  the path-connected components of $\cS_\cP$   define a foliation denoted by $\FP$.   That is, 
  $n$-dimensional weak solenoids   are examples of $n$-dimensional  matchbox manifolds. The properties of weak solenoids have been studied   in the works by Schori  \cite{Schori1966}, Rogers and Tollefson \cite{RogersJT1970,RT1971a,RT1971b,RT1972a},    Fokkink and Oversteegen in \cite{FO2002}, and by the authors.

  Given a basepoint $x \in \cS_{\cP}$   for  each $\ell \geq 0$ we then obtain   basepoints  $x_{\ell} = \Pi_{\ell}(x) \in M_{\ell}$. Denote by   $G_0 =  \pi_1(M_{0}, x_{0})$ the fundamental group of $M_0$ with basepoint $x_0$, and for $\ell > 0$  let  
\begin{equation}\label{eq-imahes}
G_{\ell} = {\rm image}\left\{  p^0_{\ell}  \colon  \pi_1(M_{\ell}, x_{\ell})\longrightarrow G_{0}\right\}  
\end{equation}
 where by a small abuse of notation, we let $p^0_{\ell}$ also denote the induced map on fundamental groups, and also suppress the dependence on the choice of basepoint $x$ in the notation $G_{\ell}$.
 Note that each $G_{\ell}$ is a subgroup of finite index in $G_0$.
  In this way, associated to the presentation $\cP$ and basepoint $x \in \fX_0$,  we obtain a descending chain of subgroups of finite index
  \begin{equation}\label{eq-descendingchain}
\cG_x \equiv \{  G_{0} \supset G_{1} \supset G_{2} \supset \cdots \supset G_{\ell} \supset \cdots   \} ~ .
\end{equation}
 
   Each quotient  $X_{\ell} = G_{0}/G_{\ell}$ is a finite set equipped with a left $G_0$-action, and there are surjections $X_{\ell +1} \to X_{\ell}$ which commute with the action of $G_0$.  The inverse limit,  
\begin{equation}\label{eq-Galoisfiber}
X(\cG) = \varprojlim  ~ \{ p_{\ell +1} \colon X_{\ell +1} \to X_{\ell} \mid \ell \geq 0\} = \{(eG_0, g_1G_1, \ldots) \mid g_{\ell} G_{\ell} = g_{\ell+1}G_{\ell}\} ~ \subset \prod_{\ell \geq 0} ~ X_{\ell}  
\end{equation}
is then a totally-disconnected, compact, perfect set for the topology induced from the Tychonoff topology on the product, so is   a Cantor set. 
Note that there is a canonical basepoint   $(eG_{\ell}) \in X(\cG)$ where $e \in G_0$ is the identity element of the group.
The   fundamental group $G_0$ acts on the left on  $X(\cG)$ via the coordinate-wise multiplication on the product in \eqref{eq-Galoisfiber}. 
The left action of $G_0$ on each quotient space $X_{\ell}$ is transitive, so each orbit of $G_0$  acting on $X(\cG)$ is dense.
We denote this   minimal Cantor action by $\Phi_0 \colon G_0 \times X(\cG) \to X(\cG)$, or by $(X(\cG) , G_0 , \Phi_0)$.

The space $X(\cG)$ admits a natural identification with the fiber $\fX_0 \equiv \Pi_0^{-1}(x_0)$. The Cantor action   $(X(\cG) , G_0 , \Phi_0)$ is then conjugate  with the monodromy action of $G_0$ on  $\fX_0$ induced by the leaves of the foliation $\FP$. See \cite{DHL2016c} for more details on this identification, and also the dependence on the choice of the basepoint  $x \in  \fX_0$.

  A group chain $\cG$ as in \eqref{eq-descendingchain} determines a presentation $\cP$ with pointed base manifold $(M_0, x_0)$, and thus a weak solenoid $\cS_{\cP}$. We next consider two equivalence relations on group chains, which  were introduced by Rogers and Tollefson in \cite{RT1971b} and 
  Fokkink and Oversteegen in \cite{FO2002}.

Let  $\fG$ denote  the collection of all   subgroup chains in $G_0$.
The first notion of  equivalence  for elements of $\fG$  corresponds to the standard notion of intertwined chains:
\begin{defn} (Rogers-Tollefson, \cite{RT1971b}) \label{defn-greq}
Group chains $\cG = \{G_{\ell}\}_{\ell \geq 0}$ and $\cH = \{H_{\ell}\}_{\ell \geq 0}$   are \emph{equivalent} if   $G_0=H_0$ and there is a group chain $\cK = \{K_{\ell}\}_{\ell \geq 0}$ with $K_0 = G_0$, and infinite subsequences $\{G_{\ell_k}\}_{k \geq 0}$ and $\{H_{j_k}\}_{k \geq 0}$ such that $K_{2k} = G_{\ell_k}$ and $K_{2k+1} = H_{j_k}$ for $k \geq 0$.
\end{defn}

 The second notion of equivalence   is more subtle, as it corresponds to conjugating a group chain by an element of the profinite core group $C_{\infty}$ as discussed in \cite{DHL2016a}. 
\begin{defn} (Fokkink-Oversteegen, \cite{FO2002}) \label{defn-conjequiv}
Group chains $\cG = \{G_{\ell}\}_{\ell \geq 0}$ and $\cH = \{H_{\ell}\}_{\ell \geq 0}$ in $\fG$ are \emph{conjugate equivalent}    if   
$G_0=H_0$, and:
\begin{enumerate}
\item there exists a sequence $(g_{\ell})_{\ell \geq 0}$ where each $g_{\ell} \in G_0$;
\item the compatibility condition   $g_{\ell}G_{\ell} = g_{\ell +1} G_{\ell}$ for all $\ell\geq 0$ is satisfied;
\item the group chains $\cG^{(g_{\ell})} = \{g_{\ell} G_{\ell} g_{\ell}^{-1}\}_{\ell \geq 0}$ and $\cH = \{H_{\ell}\}_{\ell \geq 0}$ are equivalent.
\end{enumerate}
\end{defn}

The dynamical meaning of the equivalences in Definitions~\ref{defn-greq} and \ref{defn-conjequiv} is given by the following theorem, which follows from results  in \cite{FO2002}; see also \cite{DHL2016a}.  

\begin{thm}\label{equiv-rel-11}
Let $\cG = \{G_{\ell}\}_{\ell \geq 0}$ and $\cH = \{H_{\ell}\}_{\ell \geq 0}$ be group chains  with $H_0 = G_0$, and let 
\begin{eqnarray*}
X(\cG) & = &  \varprojlim \, \{G_0/G_{\ell+1} \to G_0/G_{\ell}\}  \ , \\
X(\cH) & = &  \varprojlim \, \{H_0/H_{\ell+1} \to H_0/H_{\ell}\}  \ .
\end{eqnarray*}
Then:
\begin{enumerate}
\item \label{er-item1} The group chains $\cG$ and $\cH$ are \emph{equivalent} if and only if there exists a homeomorphism $\tau \colon  X(\cG) \to X(\cH)$ equivariant with respect to the $G_0$-actions on $X(\cG)$ and $X(\cH)$, and such that $\phi(eG_{\ell}) = (eH_{\ell})$.
\item The group chains $\cG$ and $\cH$ are \emph{conjugate equivalent} if and only if there exists a homeomorphism $\tau \colon  X(\cG) \to X(\cH)$ equivariant with respect to the $G_0$-actions on $X(\cG)$ and $X(\cH)$.
\end{enumerate}
\end{thm}

That is, an equivalence of two group chains corresponds to the existence of a \emph{basepoint-preserving} equivariant homeomorphism between their inverse limit systems, while a conjugate equivalence of two group chains corresponds to the existence of a equivariant conjugacy between their inverse limit systems, which need not preserve the basepoint. 

\begin{remark}\label{rmk-kernels}
{\rm 
The \emph{kernel} of a group chain  $\cG = \{G_{\ell}\}_{\ell \geq 0}$ is the subgroup of $G_0$ given by
\begin{equation}\label{eq-kernel}
K(\cG) = \bigcap_{\ell \geq 0} ~ G_{\ell} \ . 
\end{equation}
It follows immediately   from the definitions that if the group chains $\cG = \{G_{\ell}\}_{\ell \geq 0}$ and $\cH = \{H_{\ell}\}_{\ell \geq 0}$ with $G_0=H_0$ are  equivalent, then $K(\cG) = K(\cH) \subset G_0$.
If the chains $\cG$ and $\cH$ are only conjugate equivalent, then the kernels need not be equal. Thus, the kernel $K(\cG)$ of a group chain $\cG$ is an invariant of equivalence, but is not necessarily invariant by conjugate equivalence.

The kernel $K(\cG)$ has an interpretation in terms of the topology of the leaves of the foliation $\FP$ of the weak solenoid $\cS_{\cP}$.
Choose a basepoint  $x \in  \Pi_0^{-1}(x_0)$ and let $L_x \subset \cS_{\cP}$ be the leaf containing $x$.
The restriction of the bundle projection  $\Pi_0|_{L_x} \colon  L_x \to M_0$  is a covering map. 
Let  $\wtzM$ be the universal cover of $M_0$. Then by   standard arguments of covering space theory (see also McCord \cite{McCord1965}) there is a homeomorphism
$\ds  \wtzM/ K(\cG) \to L_x$. 
  In particular, if $L_x$ is   simply connected, then $K(\cG) = \{e\}$.
}
 \end{remark}

A standard fact about weak solenoids, is that their homeomorphism type remains unchanged if a finite number of the initial maps in a presentation $\cP$ as in \eqref{eq-presentation} are deleted. See \cite{CHL2018a,DHL2016c,HL2018} for   discussions of this property and its significance. 
This fact suggests   a third type of equivalence relation on group chains. For a group chain $\cG= \{G_{\ell}\}_{\ell \geq 0}$  and   $k \geq 0$, define the ``truncated group chain'' $\cG^{(k)}$ by setting where $G^{(k)}_{\ell} = G_{\ell + k}$, then 
  \begin{equation}\label{eq-truncatedchain1}
\cG^{(k)} \colon  G^{(k)}_{0} \supset G^{(k)}_{1} \supset G^{(k)}_{2} \supset \cdots \supset G^{(k)}_{\ell} \supset \cdots   \, .
\end{equation}

\begin{defn}   \label{defn-return}
Group chains $\cG= \{G_{\ell}\}_{\ell \geq 0}$ and $\cH = \{H_{\ell}\}_{\ell \geq 0}$   are \emph{return equivalent} if there exists integers $k, m \geq 0$ such that the  truncated group chains $\cG^{(k)}$  and $\cH^{(m)}$ are   equivalent in the sense of Definition~\ref{defn-greq}. We say they are \emph{conjugate return equivalent} if  there exists integers $k, m \geq 0$ such that the    truncated group chains $\cG^{(k)}$  and $\cH^{(m)}$ are   equivalent in the sense of Definition~\ref{defn-conjequiv}.
\end{defn}

Recall that for the Tychonoff topology on $X(\cG)$, we have a descending chain $\{U_k(\cG) \mid k =0,1,2, \ldots\}$ of clopen neighborhoods of the basepoint $(eG_{\ell})$, where  
\begin{equation}\label{eq-truncatedchain2}
U_k(\cG)  =    \varprojlim \, ~ \{ p_{\ell +1} \colon X_{\ell +1} \to X_{\ell} \mid \ell \geq k\} 
\end{equation}
The spaces $U_k(\cG)$ and $X(\cG^{(k)})$ are canonically isomorphic, and the action of $G_k$ on $X(\cG^{(k)})$ is conjugate to the restricted action of $G_k$ on $U_k(\cG)$. Thus, in the terminology of \cite{HL2018} we have:
\begin{prop}\label{prop-re}
Given group chains $\cG= \{G_{\ell}\}_{\ell \geq 0}$ and $\cH = \{H_{\ell}\}_{\ell \geq 0}$, then:
\begin{enumerate}
\item $\cG$ and $\cH$    are  return equivalent if an only if  
the minimal Cantor actions   $(X(\cG) , G_0 , \Phi_0)$ and $(X(\cH) , H_0 , \Psi_0)$ are return equivalent by a basepoint preserving local homeomorphism. 
\item $\cG$ and $\cH$   are  conjugate return equivalent if and only if  
the minimal Cantor actions   $(X(\cG) , G_0 , \Phi_0)$ and $(X(\cH) , H_0 , \Psi_0)$ are return equivalent.
\end{enumerate}
 \end{prop}
 
Finally, we recall a basic result  from \cite{CHL2018a} that was mentioned above as the motivation for introducing the equivalence relation in Definition~\ref{defn-return}:
\begin{thm}\label{thm-re}
Let $\cS_{\cP}$ and $\cS_{\cP'}$ be weak solenoids associated to presentations $\cP$ and $\cP'$ with a common base manifolds, $M_0 = M_0'$.
Let $\cG$ be the group chain associated to $\cP$, and $\cG'$ be the group chain associated to $\cP'$.  If  $\cS_{\cP}$ and $\cS_{\cP'}$ are homeomorphic as continua, then the     group chains $\cG$ and $\cG'$ are conjugate   return equivalent. 
\end{thm}
There is a converse to this result, which requires additional assumptions on the base manifold $M_0$ and the geometry of the foliations, as given in \cite[Theorem 1.5]{CHL2018a}.
The conclusion of Theorem~\ref{thm-re} is one of the   motivations for the study of   return equivalence of group chains, with the goal of obtaining computable  invariants. The \emph{asymptotic discriminant}, introduced in \cite{HL2017} and further studied in \cite{HL2018}, is an invariant of  conjugate return equivalence. The \emph{growth of orbits} as discussed in \cite{DHL2016b} is another such invariant.

\section{Induced maps between continua}\label{sec-induced}

We recall some results from the literature, concerning the properties of induced maps between continua as defined by maps between their presentations as in \eqref{eq-presentation-top}. The notion of equivalence is a fundamental aspect of the theory of shape equivalence, as developed for example in \cite{MardesicSegal1982}. In order to apply these results to the group chains associated to   presentations, care must be taken with respect to basepoints. For this reason, we take our definitions from the original works by Bousfield and Kan \cite{BousfieldKan1972}, Fort and McCord \cite{FM1966}, Mioduszewski \cite{Mioduszewski1963}, and Rogers and Tollefson \cite{RT1971a,RT1971b}. These results are then applied for the presentations associated to a matchbox manifold which is $M$-like.

Assume there are given presentations
 \begin{eqnarray} 
\cP & = &  \{\, f_{\ell+1} \colon  X_{\ell+1} \to X_\ell \mid  \ell \geq 0\} \\
\cQ & = &  \{\, g_{\ell+1} \colon  Y_{\ell+1} \to Y_\ell \mid  \ell \geq 0\} ~ ,
\end{eqnarray}
where for all $\ell \geq 0$, the spaces   $X_{\ell}$ and $Y_{\ell}$ are compact   connected polyhedra, and the maps  $f_{\ell+1}$ and $g_{\ell+1}$ are continuous surjections. Let $\cS_{\cP}$ be the inverse limit   defined by $\cP$, and $\cS_{\cQ}$ the inverse limit of $\cQ$.  
A \emph{pro-map} from $\cP$ to $\cQ$ is a collection 
$\widehat{\sigma} \equiv \{ \sigma_{\ell} \colon X_{m_{\ell}} \to Y_{\ell} \mid \ell \geq 0\}$, where
\begin{enumerate}
\item  $0 \leq m_0 <m_1 < m_2 < \cdots$ is an increasing sequence, 
\item   $\sigma_{\ell} \colon X_{m_{\ell}} \to Y_{\ell}$ are continuous   maps for $\ell \geq 0$, 
\end{enumerate}
and moreover,  for each $\ell > 0$ and  $f^{m_{\ell+1}}_{m_{\ell}} = f_{m_{\ell} +1} \circ \cdots \circ f_{m_{\ell+1}}$, the following diagram  commutes
\begin{align}\label{eq-commutingdiagram}
\xymatrix{
 X_{m_{\ell}}\ar[d]^{\sigma_{\ell}} 
& \hspace{2mm} X_{m_{\ell+1}}  \ar[d]^{\sigma_{\ell+1} ~ .}  \ar[l]_{f^{m_{\ell+1}}_{m_{\ell}}}\\
Y_{\ell}  
& \hspace{2mm} Y_{\ell+1}  \ar[l]^{g_{\ell+1}} 
} 
\end{align}
A pro-map $\widehat{\sigma}$   determines  a   continuous map $S(\widehat{\sigma}) \colon \cS_{\cP} \to \cS_{\cQ}$  called the \emph{induced map} for   $\widehat{\sigma}$. 

Recall that $\Pi_{\ell}^{\cP} \colon \cS_{\cP}  \to X_{\ell}$ and  $\Pi_{\ell}^{\cQ} \colon \cS_{\cQ}  \to Y_{\ell}$ denote the projections on the factors of these spaces.
Then the defining property for an induced map is that for all $\ell > 0$, we have 
\begin{equation}\label{eq-commutingsquare}
\Pi_{\ell}^{\cQ} \circ S(\widehat{\sigma})= \sigma_{\ell} \circ \Pi_{m_{\ell}}^{\cP} \colon \cS_{\cP} \to Y_{\ell} ~ .
\end{equation}
A fundamental question is when is there an induced map $S(\widehat{\sigma})$ which equals $h$, or is approximately equal to $h$ in an appropriate sense? 
This problem has been analyzed in the works \cite{Mioduszewski1963,FM1966,Kaul1967,RogersJT1970,RT1971a,RT1972b,Dydak1975,Dydak1980}, and we recall their conclusions as needed for the study of $M$-like spaces.

The first result  was obtained by Mioduszewski   in \cite{Mioduszewski1963}, where he introduced a generalization of the notion of a map between resolutions, where the hypothesis    that the diagram \eqref{eq-commutingdiagram} commutes is replaced by a condition that it \emph{almost commutes}.   To be more precise, first  assume that each space $Y_{\ell}$ has a metric $d_{Y_{\ell}}$ defining its topology, for $\ell \geq 0$. Let $D(Y_{\ell}) > 0$ be the diameter of $Y_{\ell}$ for this metric, then and give the inverse limit $\cS_{\cQ}$ the    metric
\begin{equation}\label{eq-hamming}
d_{\cS_{\cQ}}((y_{\ell})_{\ell \geq 0}) , (y'_{\ell})_{\ell \geq 0})) = \sum_{\ell \geq 0} ~ \frac{1}{2^{-\ell} D(Y_{\ell}) } \cdot d_{Y_{\ell}}(y_{\ell} , y'_{\ell}) ~ .
\end{equation}
Assume that there is given a descending sequence $\widehat{\epsilon} \equiv \{\e_{0} > \e_1 >  \e_2 > \cdots\}$ with $\ds \lim_{\ell \to \infty} \, \e_{\ell} = 0$. Then a sequence of maps $\widehat{\sigma}$ as above is said to be an $\widehat{\e}$-pro-map if, in place of the assumption that the diagrams \eqref{eq-commutingdiagram} are commutative, we assume that for all $\ell > 0$, 
\begin{equation}\label{eq-epromap}
d_{Y_{\ell}}(g_{\ell+1} \circ \sigma_{\ell +1}(x), \sigma_{\ell} \circ f^{m_{\ell+1}}_{m_{\ell}}(x)) < \e_{\ell} \quad {\rm for ~ all} ~ x \in X_{m_{\ell +1}} ~ .
\end{equation}

By  \cite[Theorem 2']{Mioduszewski1963}, an $\widehat{\e}$-pro-map $\widehat{\sigma}$ determines  a   continuous map $S(\widehat{\sigma}) \colon \cS_{\cP} \to \cS_{\cQ}$. Moreover, by the results of \cite[Section 3]{Mioduszewski1963}, if there is given a homeomorphism  $\Phi \colon \cS_{\cP} \to \cS_{\cQ}$ then there exists   $\widehat{\epsilon}$ as above, and  an $\widehat{\e}$-pro-map
$\widehat{\sigma}$    such that the induced map $S(\widehat{\sigma}) \colon \cS_{\cP} \to \cS_{\cQ}$ is a  homeomorphism. A natural problem is to then ask when the induced map $S(\widehat{\sigma})$ is   equal to the   map $\Phi$.
Fort and McCord   gave a partial answer to this question, in the case where the target space $\cS_{\cQ}$ is a weak solenoid:
\begin{thm}  \cite[Theorem 1]{FM1966} \label{thm-FMc}
Let  $\cP$ be a presentation whose spaces $X_{\ell}$ are compact   connected polyhedra, and let $\cQ$ be a presentation consisting of covering maps. Assume that all spaces have metrics as above. Then for any map  $\Phi \colon \cS_{\cP} \to \cS_{\cQ}$ and $\e > 0$, there exists a descending sequence $\widehat{\epsilon}$ as above, and an $\widehat{\e}$-pro-map $\widehat{\sigma}$ such that the induced map 
$S(\widehat{\sigma}) \colon \cS_{\cP} \to \cS_{\cQ}$ is $\e$-homotopic to $\Phi$. In particular, the maps $\Phi$ and $S(\widehat{\sigma})$ are $\e$-close. 
\end{thm}

 In the proof of this result, the fact that each bonding map $q_{\ell}$ is a  covering is used to lift a homotopy from the approximation at level $\ell$ to one at level $\ell+1$ (see the proof of \cite[Lemma 5]{FM1966}). 
Rogers and Tollefson gave a further refinement of Theorem~\ref{thm-FMc} in \cite{RT1971a}, characterizing the maps for which one can obtain $h = S(\widehat{\sigma})$ in terms of whether $h$ is locally fiber preserving, and in    \cite{RT1972b} they applied their method to the case where both presentations consists of covering maps.

 Next, let $\fM$ be a matchbox manifold which is $M$-like, where $M$ is a closed manifold. Then by Theorem~\ref{thm-weak}, 
   there exists a   presentation  $\cP = \{ p_{\ell +1} \colon M_{\ell +1} \to M_{\ell}  \mid  \ell \geq 0 \}$ 
 where each   bonding map  $p_{\ell +1}$ is a proper covering map, and a homeomorphism
 \begin{equation}\label{eq-ws2}
\Phi_{\cP} \colon \fM ~  \cong  ~ \cS_{\cP} \equiv \varprojlim\,   \{ p_{\ell +1} \colon M_{\ell +1} \to M_{\ell}  \mid  \ell \geq 0 \} ~.
\end{equation}
The choice of such a presentation is not unique. The result of Rogers and Tollefson \cite[Theorem 4.6]{RT1972b} relates two such presentations. 
\begin{thm}    \label{thm-RTcoverings} 
Suppose there are  two presentations as weak solenoids
 \begin{eqnarray}
\cP & = &  \{ p_{\ell +1} \colon M_{\ell +1} \to M_{\ell}  \mid  \ell \geq 0 \} \\
\cQ & = & \{ q_{\ell +1} \colon N_{\ell +1} \to N_{\ell}  \mid  \ell \geq 0 \} ~ , 
\end{eqnarray}
 a homeomorphism $\Phi  \colon    \cS_{\cP}\to     \cS_{\cQ}$ and $\e > 0$.
Then there exists a decreasing sequence $\widehat{\epsilon}$ with $\e_0 \leq \e$,  
\begin{enumerate}
\item an intertwined increasing sequence $0 \leq  j_0 < i_0 <j_1 <i_1 < j_2 < i_2 < \cdots$, 
\item covering maps $\lambda_{\ell} \colon M_{i_{\ell +1}} \to N_{j_{\ell}}$ for $\ell \geq 0$,
\item covering maps $\mu_{\ell} \colon N_{j_{\ell}} \to M_{i_{\ell}}$ for $\ell \geq 0$ ~ , 
\end{enumerate}
such that  the following diagram $\widehat{\epsilon}$-commutes:
    \begin{align} \label{eq-commutativediagram}
 \xymatrix{
 M_{i_{0}}    &   \ar[l]  M_{i_1} \ar[ld]_{\lambda_0}  &  \ar[l]  M_{j_1}     &  \ar[l]  M_{i_2} \ar[ld]_{\lambda_1}     &  \ar[l]  M_{j_2}  \cdots       \\
N_{j_{0}}    &  \ar[l]     N_{i_1}   &  \ar[l]   N_{j_1} \ar[lu]_{\mu_1}      &  \ar[l]  N_{i_2}   &  \ar[l]  N_{j_2}  \ar[lu]_{\mu_2}  \cdots       
} 
\end{align}
Moreover, for $\widehat{\lambda} = \{\lambda_{\ell} \colon M_{i_{\ell +1}} \to N_{j_{\ell}} \mid \ell \geq 0\}$, the induced map 
$S(\widehat{\lambda})$ is $\e$-homotopic to $\Phi$, and similarly  the induced map $S(\widehat{\mu})$ is $\e$-homotopic to $\Phi^{-1}$.  
\end{thm}

We next recall an observation from Fokkink and Oversteegen in \cite{FO2002}, based on \cite[Theorem 2]{RT1971a}. Let $\cP$ and $\cQ$ be presentations for weak solenoids and  $\Phi  \colon    \cS_{\cP} \to     \cS_{\cQ}$ a homeomorphism, as in Theorem~\ref{thm-RTcoverings}. 
 Choose a basepoint  $x \in \cS_{\cP}$ and set $y = \Phi(x) \in \cS_{\cQ}$.

For $\ell \geq 0$, let $\Pi^{\cP}_{\ell} \colon \cS_{\cP} \to M_{\ell}$ be the projection map onto the factor space $M_{\ell}$, and set 
$x_{\ell} = \Pi^{\cP}_{\ell}(x) \in M_{\ell}$. Then $X(\cP)_{\ell} = (\Pi^{\cP}_{\ell})^{-1}(x_{\ell}) \subset \cS_{\cP}$.
Similarly,   let     $\Pi^{\cQ}_{\ell} \colon \cS_{\cQ} \to N_{\ell}$ be the projection map onto the factor space $N_{\ell}$,   and set 
$y_{\ell} = \Pi^{\cQ}_{\ell}(y) \in N_{\ell}$. Then $X(\cQ)_{\ell} = (\Pi^{\cQ}_{\ell})^{-1}(y_{\ell}) \subset \cS_{\cQ}$.

\begin{remark}\cite[Lemma 21]{FO2002} \label{rmk-FO}
{\rm 
In the conclusion of Theorem~\ref{thm-RTcoverings}, for $j_0 \geq 0$ and $i_1 > j_0$ sufficiently large, we can choose $\widehat{\lambda}$  so that the induced map satisfies 
$S(\widehat{\lambda})(X(\cP)_{i_1})  =  X(\cQ)_{j_0}$.
}
\end{remark}

  Fokkink and Oversteegen  use this remark  in \cite{FO2002}  to show that, in essence,    associated to a weak solenoid $\fM$, there are  well-defined pro-fundamental groups. This shown in Theorem~\ref{thm-solprogroups} below. 
    
Next, we recall the following sharper version of Theorem~\ref{thm-Freudenthal}, due to   Marde\v{s}i\'{c} and  Segal:
  \begin{thm}   \cite[Theorem~1*]{MardesicSegal1963} \label{thm-MardSeg}
Let    $\fM$  be  a continuum which is $M$-like for some closed connected manifold $M$. Then there exists a presentation 
$\cQ   =   \{ q_{\ell +1} \colon M \to M  \mid  \ell \geq 0 \}$, where each map $q_{\ell +1}$ is a continuous surjection,   and    a homeomorphism
$\Phi_{\cQ} \colon \fM    \cong    \cS_{\cQ}$.
\end{thm}

For an $M$-like matchbox manifold $\fM$, there is a presentation $\cP$  as a weak solenoid and a homeomorphism $\Phi_{\cP} \colon \fM \to  \cS_{\cP}$, as in \eqref{eq-ws2}, and a presentation $\cQ$ and a homeomorphism $\Phi_{\cQ} \colon \fM \to  \cS_{\cQ}$, as in Theorem~\ref{thm-MardSeg}. Set $\Phi =  \cS_{\cQ} \circ \cS_{\cP}^{-1} \colon \cS_{\cP} \to \cS_{\cQ}$.
Combining the results of this section, we then have:
\begin{thm} \label{thmRT} 
Suppose there are   presentations 
 \begin{eqnarray}
\cP & = &  \{ p_{\ell +1} \colon M_{\ell +1} \to M_{\ell}  \mid  \ell \geq 0 \} \label{eq-presentationWS}\\
\cQ & = & \{ q_{\ell +1} \colon M \to M \mid  \ell \geq 0 \} \label{eq-presentationML}
\end{eqnarray}
where for $\ell \geq 0$, each $p_{\ell +1}$ is a covering map,   each $q_{\ell +1}$ is a continuous surjection, and there is a homeomorphism  $\Phi \colon \cS_{\cP} \to     \cS_{\cQ}$.
Assume that the presentations have metrics as in \eqref{eq-hamming}. Then for $\e > 0$ there exists a descending sequence $\widehat{\epsilon}$ with $0 < \e_0 < \e$ and: 
\begin{enumerate}
\item intertwined increasing sequence $0 \leq  j_0 < i_0 <j_1 <i_1 < j_2 < i_2 < \cdots$, 
\item continuous surjections $\lambda_{\ell} \colon M_{i_{\ell +1}} \to M$ for $\ell \geq 0$,
\item continuous surjections $\mu_{\ell} \colon M \to M_{i_{\ell}}$ for $\ell \geq 0$,
\end{enumerate}
such that  the  following diagram $\widehat{\epsilon}$-commutes:
    \begin{align} \label{eq-commutativediagram2}
 \xymatrix{
 M_{i_{0}}    &   \ar[l]  M_{i_1} \ar[ld]_{\lambda_0}  &  \ar[l]  M_{j_1}     &  \ar[l]  M_{i_2} \ar[ld]_{\lambda_1}     &  \ar[l]  M_{j_2}  \cdots       \\
M    &  \ar[l]     M   &  \ar[l]   M \ar[lu]_{\mu_1}      &  \ar[l]  M   &  \ar[l]  M  \ar[lu]_{\mu_2}  \cdots       
} 
\end{align}

Moreover, for $\widehat{\lambda} = \{\lambda_{\ell} \colon M_{i_{\ell +1}} \to M \mid \ell \geq 0\}$, the induced map 
$S(\widehat{\lambda}) \colon \cS_{\cP} \to \cS_{\cQ}$ is $\e$-close to $\Phi$. 

For   $\widehat{\mu} = \{\mu_{\ell} \colon M \to M_{i_{\ell}} \mid \ell \geq 0\}$,   the induced map 
$S(\widehat{\mu}) \colon \cS_{\cQ} \to \cS_{\cP}$ is $\e$-homotopic to $\Phi^{-1}$,  and in particular the maps $\Phi^{-1}$ and $S(\widehat{\mu})$ are $\e$-close.  
\end{thm}

\begin{remark}\label{rmk-compare}
{\rm 
In  Section~\ref{sec-progroupsWS}, we compare in Theorem~\ref{thm-solprogroups} the pro-groups obtained from the fundamental groups for the spaces in the top and bottom rows of \eqref{eq-commutativediagram} in Theorem~\ref{thm-RTcoverings}, and the fact that we are given that the squares formed by the maps in $\widehat{\lambda}$ and $\widehat{\mu}$ are $\widehat{\e}$-homotopic, and that the maps in   $\widehat{\lambda}$ and $\widehat{\mu}$ are coverings, suffices to imply that these pro-groups are isomorphic. 

Then in Section~\ref{sec-progroupsMlike}, we compare in Theorem~\ref{thm-proM} the pro-groups obtained from the fundamental groups for the spaces in the top and bottom rows of \eqref{eq-commutativediagram2} in Theorem~\ref{thmRT}. However, in this case, it is no longer given that the maps in   $\widehat{\lambda}$ and $\widehat{\mu}$ are coverings. We note that   the proof of Theorem~\ref{thmRT} by Marde\v{s}i\'{c} and Segal shows that there are simplicial spaces  $N_{\ell}$, with each $N_{\ell}$ homeomorphic to $M$, such that the maps 
$q_{\ell +1} \colon M \to M $ for $\ell \geq 0$ are induced by simplicial maps 
$q'_{\ell +1} \colon N_{\ell+1} \to N_{\ell}$. Thus, we can use the results of 
James W. Rogers in \cite{RogersJW1973} to conclude that the   squares formed by the maps in $\widehat{\lambda}$ and $\widehat{\mu}$ are in fact $\widehat{\e}$-homotopic. However, it is not given that the maps   in $\widehat{\lambda}$ and $\widehat{\mu}$ are locally onto, as they need not be covering maps, so that one does not have the injectivity property on fundamental groups that covering maps enjoy. Thus, an additional argument is needed to compare the pro-groups obtained from the fundamental groups for the spaces in the top and bottom rows of \eqref{eq-commutativediagram2}. This additional fact is given by  Proposition~\ref{prop-leafonto}, which shows the local onto property for maps restricted to    leaves, which is then used to establish an $\e$-homotopy lifting property  that is sufficient to compare these pro-groups by the path-lifting property.
}
\end{remark}

\section{Pro-groups and group chains}\label{sec-progroups}
 
 In this section, we recall the general definition of a pro-group, and equivalence between pro-groups, following the development in \cite[Chapter II, Section 2]{MardesicSegal1982}. We then compare this notion with the group chains introduced in the previous section.

  \begin{defn}\label{defn-progroup}
  A pro-group is a collection  
  $$\cG(G_{\lambda}, p^{\lambda}_{\lambda'} , \Lambda) = \{ p^{\lambda}_{\lambda'}  \colon G_{\lambda'} \to G_{\lambda} \mid \lambda \leq \lambda' , ~ \lambda, \lambda' \in \Lambda \}$$
  where   the $G_{\lambda}$ are groups, $p^{\lambda}_{\lambda'} $ are group homomorphisms, and $\Lambda$ is a directed  index set.
  \end{defn}
  A \emph{morphism}  $(f_{\gamma}, \phi)$ between pro-groups $\cG(G_{\lambda}, p^{\lambda}_{\lambda'}  , \Lambda)$ and $\cG(H_{\lambda}, q^{\gamma}_{\gamma'} , \Gamma) $ is given by an order-preserving map $\phi \colon \Gamma \to \Lambda$ and homomorphisms $f_{\gamma} \colon   G_{\phi(\gamma)} \to H_{\gamma}$ 
  such that, whenever $\gamma < \gamma'$, then there is $\lambda \in \Lambda$, $\lambda \geq \phi(\gamma)$ and $\lambda \geq  \phi(\gamma')$,   for which 
 $f_{\gamma} \circ p^{\phi(\gamma)}_{\lambda}  = q^{\gamma}_{\gamma'} \circ f_{\gamma'} \circ p^{\phi(\gamma')}_{\lambda}$. That is, the following   commutes:
    \begin{align} \label{eq-promorphism}
 \xymatrix{
 G_{\phi(\gamma)} \ar[d]_{f_{\gamma}}   &  & \ar[ll]_{p^{\phi(\gamma)}_{\lambda}}  G_{\lambda} \ar[rr]^{p^{\phi(\gamma')}_{\lambda}}   &   & G_{\phi(\gamma')}    \ar[d]^{f_{\gamma'}}    \\
H_{\gamma}   &  &  & & \ar[llll]_{q^{\gamma}_{\gamma'}} H_{\gamma'}     
} 
\end{align}
The diagram \eqref{eq-promorphism} is called an \emph{equalizer} between the maps $f_{\gamma}$ and $f_{\gamma'}$.

 The composition of morphisms of pro-groups is again a morphism, and the operation is associative. The details of the proofs of these statements, along with the notion of equivalence between morphisms, can be found in  \cite[Chapter I, Section 1.1]{MardesicSegal1982}. Also,  the notion of monomorphism,  epimorphism and isomorphism between pro-groups is defined there. We recall a reformulation of these notions using the results in \cite[Chapter I, Sections 1.1, 1.2]{MardesicSegal1982}. 

Given a morphism $(f_{\gamma}, \phi)$, a pair $(\lambda, \gamma) \in \Lambda \times \Gamma$ is said to be \emph{admissible} if $\lambda \geq \phi(\gamma)$.
For an admissible pair $(\lambda, \gamma)$ define the homomorphism
\begin{equation}\label{eq-compos}
f^{\gamma}_{\lambda} = f_{\gamma} \circ p^{\phi(\gamma)}_{\lambda} ~ .
\end{equation}

  \begin{defn}\label{defn-promono}
  \cite[Chapter II, Section 2.1, Theorem 2]{MardesicSegal1982}
  A morphism $(f_{\gamma}, \phi)$ between pro-groups $\cG(G_{\lambda}, p^{\lambda}_{\lambda'} , \Lambda)$ and $\cG(H_{\lambda}, q^{\gamma}_{\gamma'} , \Gamma) $ is a \emph{monomorphism} if for every admissible pair $(\lambda, \gamma)$ for $(f_{\gamma}, \phi)$, there is an admissible pair  $(\lambda', \gamma')$ for $(f_{\gamma}, \phi)$ with $\lambda' \geq \lambda$ and $\gamma' \geq \gamma$, such that 
\begin{equation}\label{eq-promono}
{\rm Ker}\{f^{\gamma'}_{\lambda'} \colon G_{\lambda'} \to H_{\gamma'} \} \subset {\rm Ker}\{p^{\lambda}_{\lambda'} \colon G_{\lambda'} \to G_{\lambda}\} ~ .
\end{equation}
\end{defn}

\medskip

  \begin{defn}\label{defn-proepi}
  \cite[Chapter II, Section 2.1, Theorem 4]{MardesicSegal1982}
  A morphism $(f_{\gamma}, \phi)$ between pro-groups $\cG(G_{\lambda}, p^{\lambda}_{\lambda'} , \Lambda)$ and $\cG(H_{\lambda}, q^{\gamma}_{\gamma'} , \Gamma) $ is an \emph{epimorphism} if   for every admissible pair   $(\lambda, \gamma)$ for $(f_{\gamma}, \phi)$, there is an admissible pair  $(\lambda', \gamma')$ for $(f_{\gamma}, \phi)$ with $\lambda' \geq \lambda$ and $\gamma' \geq \gamma$, such that 
\begin{equation}\label{eq-proepi}
{\rm Image}\{q^{\gamma}_{\gamma'} \colon H_{\gamma'} \to H_{\gamma} \} \subset {\rm Image}\{f^{\gamma}_{\lambda} \colon G_{\lambda} \to H_{\gamma}\} ~ .
\end{equation}
\end{defn}

\medskip
The following characterization of isomorphism of pro-groups is then equivalent to Morita's Lemma:
  \begin{defn}\label{defn-proiso}
  \cite[Chapter II, Section 2.2, Theorem 6]{MardesicSegal1982}
  A morphism $(f_{\gamma}, \phi)$ between pro-groups $\cG(G_{\lambda}, p_{\lambda, \lambda'} , \Lambda)$ and $\cG(H_{\lambda}, q_{\gamma, \gamma'} , \Gamma) $ is an \emph{isomorphism} if and only if  $(f_{\gamma}, \phi)$ is   a monomorphism and an epimorphism.
  \end{defn}

 For our purposes, the ordered sets $\Lambda$ and $\Gamma$ above will always be assumed to be a subset of the non-negative integers $\mN$ with the natural order.

 Finally, we compare the notions of group chains and the pro-groups they determine.
 Let $\cG = \{G_{\ell}\}$ and $\cH = \{H_{\ell}\}$ be two group chains, which determine pro-groups
 \begin{eqnarray*}
\cG(G_{\ell}, p^{\ell}_{\ell'}) & = & \{ p^{\ell}_{\ell'}  \colon G_{\ell'} \to G_{\ell} \mid \ell' \geq \ell \geq 0, ~ \ell, \ell' \in \mN \} \\
\cG(H_{\ell}, q^{\ell}_{\ell'}) & = & \{ q^{\ell}_{\ell'}  \colon H_{\ell'} \to H_{\ell} \mid \ell' \geq \ell \geq 0, ~ \ell, \ell' \in \mN \} ~ .
\end{eqnarray*}
where the maps $p^{\ell}_{\ell'}$ and $q^{\ell}_{\ell'}$ are given by the inclusions of subgroups. 

Suppose the group chains $\cG$ and $\cH$ are equivalent in the sense of Definition~\ref{defn-greq}. Then there are infinite increasing sequences $\{\ell_k\}$ and $\{j_k\}$, and inclusion maps 
$$\lambda_k \colon K_{2k+2} = G_{\ell_{k+1}} \subset K_{2k+1} = H_{j_k} \quad  {\rm  and} \quad \mu_k \colon K_{2k+1} = H_{j_k} \subset  K_{2k} = G_{\ell_k} ~ .$$
Introduce the group chains $\cA = \{A_k \equiv G_{\ell_k} \mid k \geq 1\}$ and $\cB = \{B_k \equiv H_{j_k} \mid k \geq 1\}$. Then we obtain  a pro-group $\cG(A_k, \iota^k_{k'})$, where $\iota^k_{k'} \colon A_{k'} \subset A_k$ is the inclusion map.  It is immediate that $\cG(G_{\ell}, p^{\ell}_{\ell'})$ and $\cG(A_k, \iota^k_{k'})$ are isomorphic   pro-groups.
Likewise, form the  pro-group $\cG(B_k, \xi^k_{k'})$ where  $\xi^k_{k'} \colon B_{k'} \subset B_k$ is the inclusion map.  It is immediate that $\cG(H_{\ell}, q^{\ell}_{\ell'})$ and $\cG(B_k, \xi^k_{k'})$ are isomorphic   pro-groups.

For   $k \geq 1$, set $\phi(k) = k+1$, and   let $f_{k} = \lambda_k \colon A_{k+1} \to B_k$. Then the map of pro-groups $(f_k, \phi)$ is a monomorphism, as each $\lambda_k$ is a monomorphism. This map is also an epimorphism, which follows using the maps  $\{\mu_k\}$.
It follows that the pro-groups $\cG(G_{\ell}, p^{\ell}_{\ell'}) $ and $\cG(H_{\ell}, q^{\ell}_{\ell'})$ are isomorphic.

If the group chains $\cG$ and $\cH$ are equivalent in the sense of Definition~\ref{defn-conjequiv}, then there exists
  a sequence $(g_{\ell})_{\ell \geq 0}$ where each $g_{\ell} \in G_0$ and    $g_{\ell}G_{\ell} = g_{\ell +1} G_{\ell}$ for all $\ell\geq 0$, so that the groups chains $\cG^{(g_{\ell})} = \{g_{\ell} G_{\ell} g_{\ell}^{-1}\}_{\ell \geq 0}$  and $\cH$ are equivalent. Define 
\begin{equation}\label{eq-chainstopro}
\lambda_k(g) = g_{\ell_{k+1}} g g_{\ell_{k+1}}^{-1} ~ {\rm for}~ g \in G_{\ell_{k+1}} \quad , \quad 
  \mu_k(h) = g_{j_k}^{-1} h g_{j_k} ~ {\rm for} ~ h \in H_{j_k} ~ ,
\end{equation}
 then as before, this yields   an isomorphism between the pro-groups $\cG(G_{\ell}, p^{\ell}_{\ell'}) $ and $\cG(H_{\ell}, q^{\ell}_{\ell'})$.
   
  Next, suppose that we have an isomorphism $(f_{\gamma}, \phi)$ between pro-groups $\cG(G_{\ell}, p^{\ell}_{\ell'})$ and $\cG(H_{\ell}, q^{\ell}_{\ell'})$.  
Then we have monomorphisms $f_{\ell} \colon G_{\phi(\ell)} \to H_{\ell}$.
 Define $H'_{\ell} = f_{\ell}(G_{\phi(\ell)})$, then by   the commutative properties in   diagram \eqref{eq-promorphism} and the isomorphism conditions \eqref{eq-promono} and \eqref{eq-proepi}, we obtain a group chain 
 \begin{equation}
H_0 \supset H_1 \supset H'_1 \supset H_2 \supset H'_2 \supset \cdots  ~ .
\end{equation}
It follows that the group chains $\cH = \{H_{\ell}\}$ and $\cH' = \{H'_{\ell}\}$ are equivalent in the sense of Definition~\ref{defn-greq}. 

Now suppose that the pro-groups are obtained from group chains $\cG$ and $\cH$ as above. For simplicity, assume that $G_0 = H_0$. Assume that the pro-groups are isomorphic, then by the above, we have that $\cH$ and $\cH'$ are equivalent group chains. 
However, this does not imply that the group chains $\cG$ and $\cH$ in $G_0$ are conjugate equivalent, as the conditions in Definition~\ref{defn-conjequiv} require that each monomorphism $f_{\ell}$ be implemented by conjugation by some element   $g_{\ell} \in G_0$. In other words, given two isomorphic subgroups of a group, they need not be in the same conjugacy class.
If we consider the weaker notion of return equivalence in Definition~\ref{prop-re}, then the requirement becomes that there exists some $\ell_0 > 0$ such that the maps $f_{\ell}$ are implemented by a conjugacy for $\ell \geq \ell_0$ and this need not hold in general either. Thus, the notion of isomorphism between the pro-groups defined by group chains $\cG$ and $\cH$ may be weaker than the notion of return equivalence between the group chains.
 
 \begin{quest}
 Is there a natural algebraic condition on a group $G$, which implies that two group chains $\cG$ and $\cH$ in $G$ which are isomorphic as pro-groups must be return equivalent?
 \end{quest}

\section{Pro-fundamental groups of weak solenoids}\label{sec-progroupsWS}

 A presentation for a continuum $\fM$ can be viewed as a particular type of shape expansion for $\fM$, and so defines its pro-homotopy type, in the sense of  Marde\v{s}i\'{c} and  Segal 
\cite{MardesicSegal1982}, or in the foundational work by Bousfield and Kan \cite{BousfieldKan1972}. 
  One can apply various functors to this shape expansion to obtain invariants of $\fM$. Marde\v{s}i\'{c} and  Segal  show in  \cite[Chapter II, Section 3.1]{MardesicSegal1982} that the pro-homology groups {\it pro}-$H_k(\fM)$ for $k > 0$ are well-defined, and show in   \cite[Chapter II, Section 3.3]{MardesicSegal1982} that the   pro-homotopy groups {\it pro}-$\pi_k(\fM, z)$ for $k > 1$ are well-defined.
  
 The definition of   the pro-fundamental group  {\it pro}-$\pi_1(\fM, z)$  for a weak solenoid requires extra care, due to basepoint issues. If one assumes that there is a shape expansion which preserves basepoints, and that maps between expansions also preserves basepoints, then the results of  \cite[Chapter II, Section 3.3]{MardesicSegal1982} show that the pro-fundamental group is well-defined. However, the maps between presentations in Section~\ref{sec-induced} above are not assumed to preserve basepoints. The work of Fox in 
 \cite{Fox1972} introduces the notion of \emph{tropes} to circumvent these difficulties, and define covering spaces for shape approximations. This topic is addressed further in the work of 
 Marde\v{s}i\'{c} and Matijevi\'{c}   \cite{MardesicMatijevic2001}, and the work of    
 Eda,  Mandi\'{c} and Matijevi\'{c} in \cite{EMM2005} illustrates some of the subtleties of the issues involved. In any case, this approach does not suffice for the study of the equivalence of pro-fundamental groups for the inverse limits spaces we consider in this work. Instead, we develop   an approach based on the fact that the shape expansions we work with are derived from a matchbox manifold, so that we can use the geometry of the leaves of the foliation of this space to show the pro-groups are well-defined.

 Now assume that there is a solenoidal presentation $\cP$ and a homeomorphism $\Phi_{\cP} \colon \fM \to \cS_{\cP}$.
Choose a basepoint $z \in \fM$, then    set $x = \Phi_{\cP}(z)$.   
Let $\Pi^{\cP}_{\ell} \colon \cS_{\cP} \to M_{\ell}$ be the projection map onto the factor space $M_{\ell}$,   set 
$x_{\ell} = \Pi^{\cP}_{\ell}(x) \in M_{\ell}$ and let  $X(\cP)_{\ell} = (\Pi^{\cP}_{\ell})^{-1}(x_{\ell}) \subset \cS_{\cP}$.
Define the pro-group 
\begin{equation}\label{eq-weakprogroupdef}
\cG(\cP,x) = \{p_{\ell +1} \colon \pi_1(M_{\ell +1} , x_{\ell+1}) \to \pi_1(M_{\ell} , x_{\ell})   \mid  \ell \geq 0 \}
\end{equation}
where by a small abuse of notation, we let $p_{\ell +1}$ also denote the induced map on fundamental groups.
 
We then have the following result, which can be considered a formal restatement of results in \cite{FO2002}:
 \begin{thm}\label{thm-solprogroups}
Let $\fM$ be an $M$-like matchbox manifold, and let $\cP$ be a presentation for which there is a homeomorphism 
$\Phi_{\cP} \colon \fM \to \cS_{\cP}$. Choose $z \in \fM$ and set $x = \Phi_{\cP}(z)$. Then the isomorphism class of the pro-group {\it pro}-$\pi_1(\fM,z)$ derived from $\cG(\cP,x)$   is independent of the choice of presentation $\cP$ and basepoint $x$, and thus is   invariant of the homeomorphism class of $\fM$.
 \end{thm}
 \proof
 Suppose that $\cQ$ is another presentation for $\fM$ as a weak solenoid, then we are given the hypotheses  and conclusions of  Theorem~\ref{thm-RTcoverings}. We use the notation there. 
 Choose a basepoint $z \in \fM$, let $x = \Phi_{\cP}(z)$ and $y = \Phi_{\cQ}(z)$. Set $\Phi = \Phi_{\cQ}^{-1} \circ  \Phi_{\cP}$. 

For the pro-morphisms $\widehat{\lambda}$ and $\widehat{\mu}$ obtained from the maps in diagram \eqref{eq-commutativediagram} which $\widehat{\e}$-commutes, the induced map  $S(\widehat{\lambda})$ is $\e$-homotopic to $\Phi$, and similarly  the induced map $S(\widehat{\mu})$ is $\e$-homotopic to $\Phi^{-1}$.  

For $\ell \geq 0$ set $x_{\ell} = \Pi^{\cP}_{\ell}(x) \in M_{\ell}$ and $y_{\ell} = \Pi^{\cQ}_{\ell}(y) \in N_{\ell}$, then   we have:
\begin{equation}\label{eq-ehomotopic}
d_{N_{j_{\ell}}}( \lambda_{\ell}(x_{i_{\ell +1}}) ,  y_{j_{\ell}}) \leq \e_{\ell}  \quad {\rm and} \quad  
d_{M_{i_{\ell}}}( \mu_{\ell}(y_{j_{\ell}}) , x_{i_{\ell}}) \leq {\e_{\ell}} ~ .
\end{equation}
The  $\e$-homotopy from $S(\widehat{\lambda})$  to $\Phi$ defines a path from the basepoint $\lambda_{\ell}(x_{i_{\ell +1}})$ to the basepoint   $y_{j_{\ell}}$, and similarly from the basepoint $\mu_{\ell}(y_{j_{\ell}})$ to the basepoint   $x_{i_{\ell}}$. Moreover, this homotopy induces a homotopy of paths at these basepoints, so we obtain from \eqref{eq-commutativediagram} a commutative diagram:
    \begin{align} \label{eq-commutativediagrampi1}
 \xymatrix{
   &      \pi_1(M_{i_1}  , x_{i_1})  \ar[ld]_{\lambda_0}  &  \ar[l]  \pi_1(M_{j_1}  , x_{j_1})      &  \ar[l]  \pi_1(M_{i_2} , x_{i_2})   \ar[ld]_{\lambda_1}     &  \ar[l]  \pi_1(M_{j_2}  , x_{j_2})  \cdots       \\
\pi_1(N_{j_0}   , y_{j_0})   &  \ar[l]     \pi_1(N_{i_1}   , y_{i_1})  &  \ar[l]   \pi_1(N_{j_1}  , y_{j_1}) \ar[lu]_{\mu_1}      &  \ar[l]  \pi_1(N_{i_2}   , y_{i_2})  &  \ar[l]  \pi_1(N_{j_2}   , y_{j_2}) \ar[lu]_{\mu_2}  \cdots       
} 
\end{align}
where we again abuse notation, by denoting a map and its induced map on fundamental groups the same. All maps in \eqref{eq-commutativediagrampi1} are monomorphisms, as they are induced from covering maps.

We   check that the conditions of Definitions~\ref{defn-promono} and \ref{defn-proepi} are satisfied for $G_i = \pi_1(M_i,x_i)$ and $H_j = \pi_1(N_j,y_j)$. We first translate the   conditions \eqref{eq-promono} and the conditions \eqref{eq-proepi} into the nomenclature of diagram \eqref{eq-commutativediagrampi1}.
The corresponding indexing sets are $\Lambda = \{i_1 , i_2, \ldots  \}$ and $\Gamma = \{j_1 , j_2, j_3, \ldots  \}$. Recall that these indices are chosen so that $ j_{\ell} < i_{\ell} <j_{\ell + 1} <i_{\ell + 1}$ for all $\ell \geq 0$.

The maps  $f_{\lambda}$ correspond to the maps $\lambda_{\ell} = G_{i_{\ell+1}} \to H_{j_{\ell}}$ so we have $\phi(j_{\ell}) = i_{\ell+1}$.

The admissible condition $\lambda \geq \phi(\gamma)$ then becomes, for $\lambda = k$ and  $\gamma = \ell$,  that $i_{k} \geq i_{\ell +1}$ so by our choice of the sequences, that $k \geq \ell+1$.

The maps $f_{\gamma} \colon   G_{\phi(\gamma)} \to H_{\gamma}$ correspond to the maps $\lambda_{\ell} \colon G_{i_{\ell +1}} \to H_{j_{\ell}}$, and so the maps $f^{\gamma}_{\lambda} = f_{\gamma} \circ p^{\phi(\gamma)}_{\lambda}$ defined in 
\eqref{eq-compos} correspond to the maps $f^{\ell}_{k} = \lambda_{\ell} \circ p^{\phi(j_\ell)}_{k}$.
 Then the monomorphism condition   \eqref{eq-promono}   becomes
\begin{equation}\label{eq-monoaltnotation}
{\rm Ker}\{\lambda_{\ell'} \circ p^{\phi(j_{\ell'})}_{k'}   \colon G_{k'} \to H_{\ell'} \} \subset {\rm Ker}\{p^{k}_{k'} \colon G_{k'} \to G_{k}\}
\end{equation} 
  and the epimorphism condition   \eqref{eq-proepi}   becomes
\begin{equation}\label{eq-epialtnotation}
{\rm Image}\{q^{\ell}_{\ell'} \colon H_{\ell'} \to H_{\ell} \} \subset {\rm Image}\{\lambda_{\ell} \circ p^{\phi(j_{\ell})}_{k} \colon G_{k} \to H_{\ell}\} ~ .
\end{equation}
The induced maps $\lambda_{\ell}$ and $p_{\ell}$ are monomorphisms, as they are induced by covering maps, 
so condition \eqref{eq-monoaltnotation} is satisfied. 
 To show  condition \eqref{eq-epialtnotation}, we use the   commutativity of the diagram \eqref{eq-commutativediagrampi1} which implies that $q^{\ell}_{\ell'} = \lambda_{\ell} \circ p^{i_{\ell +1}}_{i_{\ell'}} \circ \mu_{\ell'}$ so that \eqref{eq-epialtnotation} is satisfied for  $k = \ell' \geq \ell+1$.

 Thus, the pro-morphisms $\widehat{\lambda}$ and $\widehat{\mu}$ induce isomorphisms between the pro-groups defined by the top and bottom rows of \eqref{eq-commutativediagrampi1}.     This shows that the pro-group $\cG(\cP,x)$ defined in \eqref{eq-weakprogroupdef} is independent of the choice of presentation $\cP$.

  It remains to show that for two choices $z,z' \in \fM$, the pro-groups $\cG(\cP,x)$ and $\cG(\cP,x')$ are  isomorphic for $x = \Phi_{\cP}(z)$ and $x' = \Phi_{\cP}(z')$. 
By Remark~\ref{rmk-FO},    we can assume   that $x' \in \fX_0$.
Then by Theorem~\ref{thm-re}, the group chains $\cG_x$ and $\cG_{x'}$ are conjugate return equivalent in the sense of Definition~\ref{defn-return}.

Then as in the discussion at the end of Section~\ref{sec-progroups}, observe that the conjugation maps in Definition~\ref{defn-conjequiv} induce isomorphisms between the subgroups 
$G_{\ell}$ and $H_{\ell}$ for $\ell > 0$, and the notion of equivalence of group chains in Definition~\ref{defn-greq} induces an isomorphism of the pro-groups $\cG(\cP,x)$ and $\cG(\cP,x')$ they determine, as was to be shown.
\endproof
 
 \begin{remark}
 {\rm 
 We note a curious aspect of the conclusion of Theorem~\ref{thm-solprogroups}. 
 The isomorphism class {\it pro}-$\pi_1(\fM,z)$ is independent of the choice of basepoint used to define the group chain $\cG^x$, for $x = \Phi_{\cP}(z)$, while by 
 Remark~\ref{rmk-kernels} the kernel of the associated group chain $K(\cG^x)$ may depend on the choice of $x$.
 Thus, there is some loss of information when we pass from a group chain $\cG^x$ to the pro-group it determines. The reason is simply that the group chain $\cG^x$ contains the information on how the groups are embedded in $G_0$, while the pro-group {\it pro}-$\pi_1(\fM,x)$ does not. 
  }
 \end{remark}

\section{Pro-fundamental groups of M-like spaces}\label{sec-progroupsMlike}

In this section we prove the following, which can be considered the main result of this work.
\begin{thm}\label{thm-proM}
 Let $\fM$ be a matchbox manifold which is  $M$-like. Then {\it pro}-$\pi_1(\fM, z)$  is pro-isomorphic to the pro-group defined by  maps $\{g_{\ell} \colon G \to G \mid \ell \geq 0\}$, where $G= \pi_1(M, y)$.
 \end{thm}
 \proof  The idea is to use Theorem~\ref{thm-MardSeg} to obtain a diagram analogous to \eqref{eq-commutativediagrampi1}, and then follow the outline of the proof of Theorem~\ref{thm-solprogroups}. This requires that we first choose basepoints that ``almost commute'', so the pro-morphisms   $\widehat{\lambda}$ and $\widehat{\mu}$ are well-defined. The main issue then is to show the maps satisfy  the conditions \eqref{eq-monoaltnotation} and \eqref{eq-epialtnotation}. It follows from Theorem~\ref{thmRT}  that the maps in $\widehat{\mu}$ induce commutative squares in the diagram \eqref{eq-commutativediagram2}.
Unfortunately,   the maps $\widehat{\lambda}$ in the diagram \eqref{eq-commutativediagram2} are only known to  $\widehat{\epsilon}$-commute, and neither collection of maps $\widehat{\lambda}$ and $\widehat{\mu}$ are known to be coverings. Thus,  another approach is required, and we use a technique that   is  analogous to the approach of Rogers and Tollefson in \cite[Section 3]{RT1972b}.

We  first  recall some basic results about the leafwise geometry of foliated spaces and matchbox manifolds, as discussed in detail for example in the works \cite{CandelConlon2000,ClarkHurder2013,CHL2018b}. 

\begin{thm}\label{thm-riemannian}
Let $\fM$ be a smooth matchbox manifold with foliation $\FfM$. Then there exists a leafwise Riemannian metric for $\FfM$, such that for each $x \in \fM$, the leaf $L_x$ inherits the structure of a complete Riemannian manifold with bounded geometry, and the Riemannian metric and its covariant derivatives depend continuously on $x$ .  
\end{thm}
Bounded geometry implies, for example, that for each $x \in \fM$, there is a leafwise exponential map
$\exp^{\F}_x \colon T_x\FfM \to L_x$ which is a surjection, and the composition with the inclusion map, $\exp^{\F}_x \colon T_x\FfM \to L_x \subset \fM$, depends continuously on $x$ in the compact-open topology on maps.

Each leaf $L \subset \fM$ has a complete path-length metric, induced from the leafwise Riemannian metric:
$$\dF(x,y) = \inf \left\{\| \gamma\| \mid \gamma \colon [0,1] \to L ~{\rm is ~ piecewise ~~ C^1}~, ~ \gamma(0) = x ~, ~ \gamma(1) = y ~, ~ \gamma(t) \in L \quad \forall ~ 0 \leq t \leq 1\right\}$$
  where $\| \gamma \|$ denotes the path-length of the piecewise $C^1$-curve $\gamma(t)$. If $x,y \in \fM$   are not on the same leaf, then set $\dF(x,y) = \infty$.
  For each $x \in \fM$ and $r > 0$, let $D_{\F}(x, r) = \{y \in L_x \mid \dF(x,y) \leq r\}$.

For each $x \in \fM$, the  {Gauss Lemma} implies that there exists $\lambda_x > 0$ such that $D_{\F}(x, \lambda_x)$ is a \emph{strongly convex} subset for the metric $\dF$. That is, for any pair of points $y,y' \in D_{\F}(x, \lambda_x)$ there is a unique shortest geodesic segment in $L_x$ joining $y$ and $y'$ and  contained in $D_{\F}(x, \lambda_x)$. This standard concept of Riemannian geometry is discussed in detail in  \cite{BC1964}, and in \cite[Chapter 3, Proposition 4.2]{doCarmo1992}. 
Then for all $0 < \lambda < \lambda_x$ the disk $D_{\F}(x, \lambda)$ is also strongly convex. 
Then   we have:
\begin{lemma}\label{lem-stronglyconvex}
There exists $\lF > 0$ such that for all $x \in \fM$, $D_{\F}(x, \lF)$ is strongly convex.
\end{lemma}

Choose  a presentation  $\cQ  =   \{ q_{\ell +1} \colon M \to M \mid  \ell \geq 0 \}$ as in Theorem~\ref{thm-MardSeg}.  The manifold $M$ is assumed to be smooth  and without boundary, and  choose a Riemannian metric $d_M$ for $M$. Let $\lambda_M > 0$ be such that  for  $u \in M$,  the geodesic ball $B_M(u, \lambda_M) \subset M$   is strongly convex.

Give $\cS_{\cQ}$  the metric   \eqref{eq-hamming}, where each factor space $Y_{\ell} = M$ is given the metric $d_M$.

Choose  a presentation  $\cP   =    \{ p_{\ell +1} \colon M_{\ell +1} \to M_{\ell}  \mid  \ell \geq 0 \}$  for which   there is a homeomorphism $\Phi_{\cP} \colon \fM    \cong    \cS_{\cP}$, 
    where each $M_\ell$ is a compact connected $n$-dimensional manifold without boundary, and  each  $p_{\ell +1 }$ is a   covering map  of finite degree greater than one.
 Choose a Riemannian metric on $M_0$ such that every geodesic ball of radius at most $\e'' > 0$ is strongly convex, then for $\ell > 0$, let $M_{\ell}$ have the Riemannian metric induced by the covering map $p^0_{\ell} \colon M_{\ell} \to M_0$ in \eqref{eq-presentation}. Give $\cS_{\cP}$  the metric as in  \eqref{eq-hamming}. 
 It then follows that for  each $b \in M_0$ and closed ball $B_{M_0}(b, \e''/2) \subset M_0$, the inverse image 
$\Pi_0^{-1}(B_{M_0}(b, \e''/2)) \subset \cS_{\cP}$ is a chart for the foliation $\FP$ on the inverse limit $\cS_{\cP}$. The proof by McCord in \cite{McCord1965}  that weak solenoids are foliated spaces used just such charts.
Also, note that each leaf $L \subset   \cS_{\cP}$ is a covering of $M_0$.  The lift of the Riemannian metric on $M_0$ to $L$ defines a     path length metric on each leaf $L$ of $\FP$. 

Moreover, the homeomorphism $\Phi_{\cP}$ preserves the path components of $\fM$ and $\cS_{\cP}$ so for each leaf $L \subset \fM$ the restriction $\Phi_{\cP} | L$ is a homeomorphism onto a leaf of $\FP$. Moreover, the compactness of $\fM$ implies that each such restriction is uniformly continuous for the leafwise metrics on $\FfM$ and $\FP$. 
Thus, we can choose $\e'' > 0$ sufficiently small so that for each $w \in \cS_{\cP}$ the leafwise ball 
$B_{\cS_{\cP}}(w, \e'') \subset L_z$ is contained in the image $\Phi(B_{\fM}(u, \lF))$, where $u = \Phi_{\cP}^{-1}(w)$.

By Theorem~\ref{thmRT}, for $\e = \min \{\e'/4, \e''/4\}$, there exists   a   sequence $\widehat{\epsilon} = \{\e \geq \e_0 > \e_1 > \cdots\}$ descending to zero, and 
\begin{enumerate}
\item an intertwined increasing sequence $0 \leq  j_0 < i_0 <j_1 <i_1 < j_2 < i_2 < \cdots$, 
\item continuous surjections $\widehat{\lambda} \equiv \{\lambda_{\ell} \colon M_{i_{\ell +1}} \to M \mid \ell \geq 0\}$,
\item continuous surjections $\widehat{\mu} \equiv \{\mu_{\ell} \colon M \to M_{i_{\ell}} \mid \ell \geq 0\}$,
\item a homeomorphism $\Phi_{\cQ} \colon \fM \to \cS_{\cQ}$,
\end{enumerate}
  so that the   diagram \eqref{eq-commutativediagram3} below $\widehat{\epsilon}$-commutes, where for the convenience of the arguments to follow we use $N_{\ell} = M$ so that the indexing is indicated:

    \begin{align} \label{eq-commutativediagram3}
 \xymatrix{
 M_{i_{0}}    &   \ar[l]  M_{i_1} \ar[ld]_{\lambda_0}  &  \ar[l]  M_{j_1}     &  \ar[l]  M_{i_2} \ar[ld]_{\lambda_1}     &  \ar[l]  M_{j_2}  \cdots       \\
N_{j_{0}}    &  \ar[l]     N_{i_1}   &  \ar[l]   N_{j_1} \ar[lu]_{\mu_1}      &  \ar[l]  N_{i_2}   &  \ar[l]  N_{j_2}  \ar[lu]_{\mu_2}  \cdots       
} 
\end{align}
 This means that for each $\ell > 0$, we have 
 \begin{eqnarray} 
d_M( \lambda_{\ell} \circ \mu_{\ell +1}(y)  , q^{j_{\ell}}_{j_{\ell +1}}(y)) & \leq &  \e_{\ell} \quad {\rm for ~ all} ~ y \in N_{j_{\ell +1}} = M \label{eq-approx1ell}\\
 d_{M_{i_{\ell}}} ( \mu_{\ell} \circ \lambda_{\ell}(x) , p^{i_{\ell}}_{i_{\ell+1}}(x))   &  \leq &   \e_{\ell} \quad {\rm for ~ all} ~ y \in M_{i_{\ell+1}} ~ . \label{eq-approx2ell}
\end{eqnarray}

Choose  $z \in \fM$, let $x = \Phi_{\cP}(z)$ and $y = \Phi_{\cQ}(z)$.
Set  $x_{\ell} = \Pi^{\cP}_{\ell}(x) \in M_{\ell}$ and $y_{\ell} = \Pi^{\cQ}_{\ell}(y) \in M$.

By \eqref{eq-approx1ell} we have that $\mu_{\ell}(y_{j_{\ell}}) \in B_{M_{i_{\ell}}}(x_{i_{\ell}}, \e_{i_{\ell}})$ which is a contractible disk. Thus, we obtain a well-defined map $\mu_{\ell} \colon \pi_1(M, y_{j_{\ell}}) \to \pi_1(M_{i_{\ell}} , x_{i_{\ell}})$ where we again use the same notation for the map on spaces and the induced map on homotopy groups. Similarly, \eqref{eq-approx2ell} implies that there  well-defined map $\lambda_{\ell} \colon \pi_1(M_{i_{\ell+1}}, x_{i_{\ell +1}}) \to \pi_1(M , y_{i_{\ell}})$. Thus we obtain the following   diagram:
\begin{align} \label{eq-commutativediagrampi12}
 \xymatrix{
   &      \pi_1(M_{i_1}  , x_{i_1})  \ar[ld]_{\lambda_0}  &  \ar[l]  \pi_1(M_{j_1}  , x_{j_1})      &  \ar[l]  \pi_1(M_{i_2} , x_{i_2})   \ar[ld]_{\lambda_1}     &  \ar[l]  \pi_1(M_{j_2}  , x_{j_2})  \cdots       \\
\pi_1(M  , y_{j_0})   &  \ar[l]     \pi_1(M  , y_{i_1})  &  \ar[l]   \pi_1(M , y_{j_1}) \ar[lu]_{\mu_1}      &  \ar[l]  \pi_1(M  , y_{i_2})  &  \ar[l]  \pi_1(M  , y_{j_2}) \ar[lu]_{\mu_2}  \cdots       
} 
\end{align}
The maps in the top row of \eqref{eq-commutativediagrampi12} are all monomorphisms as they are induced by covering maps, and the squares formed by the maps in $\widehat{\mu}$ are commutative as the   induced map 
$S(\widehat{\mu}) \colon \cS_{\cQ} \to \cS_{\cP}$ is $\e$-homotopic to $\Phi^{-1}$. We claim that the squares formed by the maps in $\widehat{\lambda}$ on fundamental groups in \eqref{eq-commutativediagrampi12}  are commutative, 
and  then show that  the properties \eqref{eq-monoaltnotation} and \eqref{eq-epialtnotation} are satisfied.

The idea of the proof of these claims is as follows. Let $L \subset \fM$ be a leaf of $\FfM$, endowed with the topology defined by the path-length metric (which is called the model $\widehat{L}$ for $L$ in \cite[Section 3]{RT1972b}). Recall that $\Phi_{\cP}(L)$ is a leaf of $\FP$ as the homeomorphism $\Phi_{\cP}$ preserves path components, and then 
the restriction $\Pi^{\cP}_{\ell} | \Phi_{\cP}(L)$ is a covering map by construction, for each $\ell > 0$. Thus,   the restriction 
$\Pi^{\cP}_{\ell} \circ \Phi_{\cP} \colon L \to M_{\ell}$ is a covering map. However,  the composition 
$\Pi^{\cQ}_{\ell} \circ \Phi_{\cQ} \colon L \to M$ need not be a covering map. What we claim is that this restriction is an approximate covering map, in that it satisfies an approximate lifting property, to be defined precisely below. 

There have been many works   studying the properties of $\e$-maps between compact manifolds, starting with those of Eilenberg \cite{Eilenberg1938} and Ganea \cite{Ganea1959}, and more recently by  
Chapman and Ferry \cite{ChapmanFerry1979}, and Ferry \cite{Ferry1979}. The idea in the following, is to adapt some of the more elementary ideas from these works to the restrictions $\Pi^{\cQ}_{\ell} \circ \Phi_{\cQ} | L$. There is an immediate difficulty, in that 
 an $\e$-map $f_{\e} \colon \fM \to M$ is assumed to be onto, but it is not given that the restriction   $\Pi^{\cQ}_{\ell} \circ \Phi_{\cQ} | L$
is locally onto. However, this property  is proved in  Proposition~\ref{Aprop-coverage} below. Recall that  $\lF > 0$ is chosen so that   for all $x \in \fM$, the leafwise closed disk   $D_{\F}(x, \lF)$ is strongly convex. Also, recall that $\lambda_M > 0$ was chosen so   that for every $u \in M$, the closed disk $D_M(u,\lambda_M) \subset M$ is strictly convex.

\begin{prop} \cite[Proposition A.1]{CHL2018b} \label{Aprop-coverage}
Let $\fM$ be a matchbox manifold with leafwise Riemannian metric on $\FfM$ and let $M$ be a closed Riemannian manifold. Then 
there exists $0 < \eF \leq \lF/4$  such that,    if $f \colon \fM \to M$ is an $\eF$-map, then there exists $\delta_1 > 0$ such that for      $x_0 \in \fM$ with $w_0 = f(x_0)$,  
\begin{equation}\label{eq-onto}
D_M(w_0,\delta_1) \subset f(D_{\F}(x_0, \lF/2 )) ~ .
\end{equation}
\end{prop}
 Thus, \eqref{eq-onto} implies that every sufficiently small disk in $M$ is contained in the image of a strictly convex disk in some leaf of $\FfM$.

\begin{prop}\label{prop-leafonto}
There exists $\ell_2 > 0$ such that for $\ell \geq \ell_2$   the map     
 $\Pi^Q_{\ell} \circ \Phi_{\cQ} \colon \fM \to M$  satisfies an approximate homotopy lifting property. That is, given an integer $k \geq 1$, constant $0 < \delta_2 < \lambda_M$ and $x \in \fM$, let $L$ be the leaf through $x$,  and given     a continuous map $\sigma \colon  [0,1]^k \to M$ with $\sigma(\vec{0}) = y_{\ell}  = \Pi^Q_{\ell} \circ \Phi_{\cQ}(x)$,    then there exists a continuous map $\widehat{\sigma} \colon  [0,1]^k \to L$ with $\widehat{\sigma}(\vec{0}) = x$, for which   
\begin{equation}
d_{M}(\Pi^Q_{\ell} \circ \Phi_{\cQ} \circ \widehat{\sigma}(\vec{v}) , \sigma(\vec{v})) \leq \delta_2 \quad {\rm for ~ all} ~ \vec{v} \in [0,1]^k~ .
\end{equation}
Moreover, the maps $\sigma$ and $\Pi^Q_{\ell} \circ \Phi_{\cQ} \circ \widehat{\sigma}$ are $\delta_2$-homotopic.
\end{prop}
\proof
Choose    $0 < \e_2 < \lF/4$, then there exists $\ell_2 > 0$ sufficiently large so that $\ell \geq \ell_2$ implies that the composition $\Pi^Q_{\ell} \circ \Phi_{\cQ}$ is an $\e_2$-map. Fix a choice of $\ell \geq \ell_2$.

Let $\delta_2 \leq \lambda_M$ be as in Proposition~\ref{Aprop-coverage}, so that for all $x_0 \in \fM$ with $w_0 = f(x_0)$,  then 
\begin{equation}
D_M(w_0,\delta_2) \subset \Pi^Q_{\ell} \circ \Phi_{\cQ}(D_{\F}(x_0, \e_2 )) ~ .
\end{equation}
Let $B(\vec{v}, \delta_3)$ denote the ball of radius $\delta_3$ in $\mR^k$ for the Euclidean metric on $\mR^k$. Then choose $\delta_3 > 0$ such that for any $\vec{v} \in [0,1]^k$, we have $\sigma(B(\vec{v}, \delta_3) \cap [0,1]^k) \subset B_M(\sigma(\vec{v}), \delta_2)$.

Choose $n_3 > 0$ so that the ``small cube''' $[0,1/n_3]^k$ has diameter less than $\delta_3$. 

Starting with the given value $\sigma(0) = y_{\ell} \in M$ we use the above choices to recursively construct the lifting $\widehat{\sigma}$, small cube by small cube for a cubical decomposition of $[0,1]^k$.
This is a standard construction, usually proved for the assumption that the target space $M$ is locally contractible, with the only nuance being the estimates on the control of the maps. We omit the details.
\endproof

   We now return to the proof of Theorem~\ref{thm-proM}.   
  We are given a solenoidal presentation $\cP$ and a homeomorphism $\Phi_{\cP} \colon \fM \to \cS_{\cP}$, and a   presentation 
$\cQ   =   \{ q_{\ell +1} \colon M \to M  \mid  \ell \geq 0 \}$, where each map $q_{\ell +1}$ is a continuous surjection,   and    a homeomorphism $\Phi_{\cQ} \colon \fM    \cong    \cS_{\cQ}$.
 Let  $\lF > 0$ as in Lemma~\ref{lem-stronglyconvex}, chosen so that for all $x \in \fM$, $D_{\F}(x, \lF)$ is strongly convex.

Let $\e_{\cP} > 0$ be such that for all $\xi \in \fM$ then $B_{\cS_{\cP}}(\Phi_{\cP}(\xi) , \e_{\cP}) \subset \Phi_{\cP}(B_{\fM}(\xi, \lF/4))$. 

Let $\e_{\cQ} > 0$ be such that for all $\xi \in \fM$ then $B_{\cS_{\cQ}}(\Phi_{\cQ}(\xi) , \e_{\cQ}) \subset \Phi_{\cQ}(B_{\fM}(\xi, \lF/4))$. 

Set $\e_{\Phi} \equiv \min \{\e_{\cP}, \e_{\cQ}, \lambda_M/4\}$, and let 
 $\widehat{\epsilon} \equiv \{\e_{0} > \e_1 >  \e_2 > \cdots\}$ with $\ds \lim_{\ell \to \infty} \, \e_{\ell} = 0$
    and $\e_0 < \e_{\Phi}$.   

We are given the   diagram \eqref{eq-commutativediagram2}, and the homeomorphism $\Phi =  \cS_{\cQ} \circ \cS_{\cP}^{-1} \colon \cS_{\cP} \to \cS_{\cQ}$, for which 
   \begin{enumerate}
\item  the maps $\widehat{\lambda} =   \{\lambda_{\ell} \colon  M_{i_{\ell+1}} \to M \mid \ell \geq 0\}$ are continuous surjections which  $\widehat{\epsilon}$-commute;
\item  the maps   $\widehat{\mu} =  \{\mu_{\ell} \colon M \to M_{i_{\ell}} \mid \ell \geq 0\}$ are covering maps which     $\widehat{\epsilon}$-homotopy commute; 
\item   $S(\widehat{\mu}) \colon \cS_{\cQ} \to \cS_{\cP}$ is $\e_{\Phi}$-homotopic to $\Phi^{-1} = \Phi =   \cS_{\cP}\circ \cS_{\cQ}^{-1}$.
\end{enumerate}
 Choose a basepoint $z \in \fM$, then    set $x = \Phi_{\cP}(z)$ and $y = \Phi_{\cQ}(z)$.   
 We   show that   the squares     in \eqref{eq-commutativediagrampi123} below are commutative, 
and   that  the properties \eqref{eq-monoaltnotation} and \eqref{eq-epialtnotation} are satisfied. 
\begin{align} \label{eq-commutativediagrampi123}
 \xymatrix{
 \cdots   &   \ar[l]   \pi_1(M_{i_{\ell+1}}  , x_{i_{\ell+1}})  \ar[ld]_{\lambda_{\ell}}  &  \ar[l]  \pi_1(M_{j_{\ell+1}}  , x_{j_{\ell+1}})      &  \ar[l]  \pi_1(M_{i_{\ell+2}} , x_{i_{\ell +2}})   \ar[ld]_{\lambda_{\ell+1}}     &  \ar[l]   \cdots       \\
\pi_1(M  , y_{j_{\ell}})   &  \ar[l]     \pi_1(M  , y_{i_{\ell+1}})  &  \ar[l]   \pi_1(M , y_{j_{\ell+1}}) \ar[lu]_{\mu_{\ell+1}}      &  \ar[l]  \pi_1(M  , y_{i_{\ell+2}})  &  \ar[l]    \cdots       
} 
\end{align}
 As noted before,    that the maps   $\widehat{\mu}$ induce commutative squares in the diagram \eqref{eq-commutativediagrampi123}.
  
 We first show that the maps   $\widehat{\lambda}$ induce commutative squares in the diagram \eqref{eq-commutativediagrampi123}.
 Let $\alpha \in \pi_1(M_{i_{\ell+2}}  , x_{i_{\ell+2}})$ be represented by a closed smooth curve 
 $$\gamma_{\alpha} \colon [0,1] \to M_{i_{\ell+2}} \quad , \quad \gamma_{\alpha}(0) = \gamma_{\alpha}(1) = x_{i_{\ell+2}} ~ .$$
 Then $\lambda_{\ell} \circ p_{i_{\ell +2}}^{i_{\ell+1}}(\alpha) \in \pi_1(M  , y_{j_{\ell}}) $ 
 is represented by the closed curve 
 $\lambda_{\ell} \circ p_{i_{\ell +2}}^{i_{\ell+1}}\circ \gamma_{\alpha} \colon [0,1] \to M$, 
 and   $q_{j_{\ell}}^{j_{\ell+1}} \circ \lambda_{\ell+1}(\alpha)$ is represented by the closed curve 
 $q_{j_{\ell}}^{j_{\ell+1}} \circ \lambda_{\ell+1} \circ \gamma_{\alpha} \colon [0,1] \to M$. We claim these two curves are homotopic.
 
 Let $\widetilde{\gamma_{\alpha}} \colon [0,1] \to L$ with  $\widetilde{\gamma_{\alpha}}(0) = z$ be the lift of $\gamma_{\alpha}$ via the covering map $\Pi_{i_{\ell +2}}^{\cP} \circ \Phi_{\cP} \colon L \to M_{i_{\ell +2}}$. 
 
 Then  $\Pi_{j_{\ell +1}}^{\cQ} \circ   \Phi_{\cQ} \circ \widetilde{\gamma_{\alpha}} \colon [0,1] \to M$. It is given that   the maps $\widehat{\lambda}$ all $\widehat{\epsilon}$-commute,  which by the definition  \eqref{eq-epromap} implies that 
 \begin{equation}
d_M( \lambda_{\ell+1} \circ \Pi_{i_{\ell +2}}^{\cP} \circ \Phi_{\cP} \circ \widetilde{\gamma_{\alpha}}(t) ~ , ~  \Pi_{j_{\ell +1}}^{\cQ} \circ   \Phi_{\cQ} \circ \widetilde{\gamma_{\alpha}}(t)) ~ \leq ~ \e_{j_{\ell+1}} ~ 0 \leq t \leq 1 ~ .
\end{equation}
As $\e_{j_{\ell+1}} \leq \lambda_M/4$ by the choice of $\e_{\Phi}$ for each $0 \leq t \leq 1$, 
there is a geodesic segment between 
$\lambda_{\ell+1} \circ \Pi_{i_{\ell +2}}^{\cP} \circ \Phi_{\cP} \circ \widetilde{\gamma_{\alpha}}(t)$ 
and 
$\Pi_{j_{\ell +1}}^{\cQ} \circ   \Phi_{\cQ} \circ \widetilde{\gamma_{\alpha}}(t)$. Thus, if we form a closed curve in $M$ by joining the endpoints of $\Pi_{j_{\ell +1}}^{\cQ} \circ   \Phi_{\cQ} \circ \widetilde{\gamma_{\alpha}}$ to $y_{j_{\ell +1}}$ by a geodesic segments in 
$B_M(y_{j_{\ell +1}}, \lambda_M/4)$, 
we obtain a curve which   represents $\lambda_{\ell+1}(\alpha)$.
A similar    argument shows that  
 \begin{equation}\label{eq-composingequality}
 \Pi_{j_{\ell }}^{\cQ} \circ   \Phi_{\cQ} \circ \widetilde{\gamma_{\alpha}} 
= q_{j_{\ell+1}}^{j_{\ell}} \circ \Pi_{j_{\ell +1}}^{\cQ} \circ   \Phi_{\cQ} \circ \widetilde{\gamma_{\alpha}} \colon [0,1] \to M
\end{equation}
represents the class $\lambda_{\ell} \circ p_{i_{\ell +2}}^{i_{\ell+1}}(\alpha) \in \pi_1(M  , y_{j_{\ell}})$, while the equality \eqref{eq-composingequality} implies that it also represents the class  $q_{j_{\ell}}^{j_{\ell+1}} \circ \lambda_{\ell+1}(\alpha)$. Thus,  they are equal as was to be shown. 
 
We next show that the induced maps  $\lambda_{\ell}$ and $\mu_{\ell}$ are monomorphisms. Let $\alpha \in \pi_1(M_{i_{\ell+1}}  , x_{i_{\ell+1}})$ and suppose that  $\lambda_{\ell}(\alpha) \in \pi_1(M  , y_{j_{\ell}})$ is trivial. We claim this implies $\alpha$ is also trivial. 
As above, let $\alpha$ be represented by a closed smooth curve 
 $$\gamma_{\alpha} \colon [0,1] \to M_{i_{\ell+1}} \quad , \quad \gamma_{\alpha}(0) = \gamma_{\alpha}(1) = x_{i_{\ell+1}} ~ ,$$
and let $\widetilde{\gamma_{\alpha}} \colon [0,1] \to L$ with  $\widetilde{\gamma_{\alpha}}(0) = z$ be the lift of $\gamma_{\alpha}$ via the covering map $\Pi_{i_{\ell +1}}^{\cP} \circ \Phi_{\cP} \colon L \to M_{i_{\ell +1}}$. Then
$\Pi_{j_{\ell}}^{\cQ} \circ   \Phi_{\cQ} \circ \widetilde{\gamma_{\alpha}} \colon [0,1] \to M$ represents $\lambda_{\ell}(\alpha)$, where as before we obtain a closed curve by adjoining short geodesic segments to the endpoints. The assumption $\lambda_{\ell}(\alpha) $ is trivial class implies that $\lambda_{\ell}\circ \alpha \colon [0,1] \to M$ extends to a homotopy to the constant map, which is given by a map $H(t,s) \colon [0,1]^2 \to M$ with $H(t,0) = \lambda_{\ell}\circ \alpha(t)$ and   $H(t,1) = H(0,s) = H(1,s) = y_{j_{\ell}}$ for all $0\leq t,s \leq 1$. Then by Proposition~\ref{prop-leafonto} there is an approximate lift $\widehat{H} \colon [0,1]^2   \to L$ such that $\Pi_{j_{\ell}}^{\cQ} \circ \Phi_{\cQ} \circ \widehat{H}$ is $\e_{j_{\ell}}$ close to $H$, and hence is $\e_{j_{\ell}}$-homotopic. Then the composition 
$\Pi_{i_{\ell+1}}^{\cP} \circ \Phi_{\cP} \circ \widehat{H} \colon [0,1]^2 \to M_{i_{\ell +1}}$ defines a homotopy of $\alpha$ to an almost constant map, hence $\alpha$ is the trivial class, as was to be shown. 

The claim that the maps $\mu_{\ell}$ induce injective maps on fundamental groups follows similarly.    
 
 The claim that the properties \eqref{eq-monoaltnotation} and \eqref{eq-epialtnotation} are satisfied now follows in exactly the same way as in the proof of Theorem~\ref{thm-solprogroups}. This completes   the proof of Theorem~\ref{thm-proM}. 
 \endproof

 As remarked previously, Marde\v{s}i\'{c} and  Segal   show in   \cite[Chapter II, Section 3.3]{MardesicSegal1982} that the   pro-homotopy groups {\it pro}-$\pi_k(\fM,z)$   are well-defined for pointed shape expansions.
The key technical fact  used in this section, Proposition~\ref{prop-leafonto},  applies equally well for the study of the   higher pro-homotopy groups {\it pro}-$\pi_k(\fM, z)$ of a weak solenoid, for $k > 1$. Instead of considering lifts of paths to the leaf $L$, one considers lifts of cubes. Then essentially the same arguments as above yield  a proof of the following result:

 \begin{thm}\label{thm-proM2}
 Let $\fM$ be a matchbox manifold which is  $M$-like. For $k > 1$,  {\it pro}-$\pi_k(\fM, z)$  is  isomorphic to the pro-group defined by  maps $\{g_{\ell} \colon \pi_k(M, y) \to \pi_k(M, y) \mid \ell \geq 0\}$. In particular, 
 {\it pro}-$\pi_k(\fM, z)$  is independent of the choice of basepoint $z \in \fM$.
 \end{thm}

\section{Proofs of main theorems}\label{sec-proofs}

In this section, we show how the conclusions and proof of Theorem~\ref{thm-proM} imply Theorems~\ref{thm-main1} and \ref{thm-main2} from the introduction. We assume that    $\fM$ be an $n$-dimensional  matchbox manifold which  is $M$-like, for a closed manifold $M$. By Theorem~\ref{thm-weak}, there exists  a presentation  $\ds \cP = \{\, f_{\ell+1} \colon  M_{\ell+1} \to M_{\ell} \mid  \ell \geq 0\}$ as a weak solenoid, 
and a homeomorphism $\Phi_{\cP} \colon \fM     \to \cS_{\cP}$.
  By Theorem~\ref{thm-MardSeg}, there exists   a presentation  $\cQ  =   \{ q_{\ell +1} \colon M \to M \mid  \ell \geq 0 \}$, where all the maps $q_{\ell +1}$ are continuous surjections, and a homeomorphism   $\Phi_{\cQ} \colon \fM     \to \cS_{\cQ}$. Choose a basepoint $z \in \fM$ and set $x = \Phi_{\cP}(z)$ and $y = \Phi_{\cQ}(z)$.
    By the proof of Theorem~\ref{thm-proM}, we can assume that all of the squares in the diagram  \eqref{eq-commutativediagrampi12} are commuting.
 
 Define the sub-presentation $\cP'$ of $\cP$ given by the maps
 \begin{equation}
\cP' = \{p_{i_{\ell+1}}^{i_{\ell}} \colon M_{i_{\ell+1}} \to M_{i_{\ell}} \mid \ell \geq 0\} ~ .
\end{equation}
Then by the basic properties of inverse limits, there is a canonical homeomorphism $\Phi' \colon \cS_{\cP'} \to \cS_{\cP}$.  

From  the proof of Theorem~\ref{thm-proM},  the maps $\lambda_{\ell}$ and $\mu_{\ell}$ in the diagram \eqref{eq-commutativediagrampi12} are monomorphisms, for all $\ell \geq 1$. Thus,  we have the extended commutative diagram of group monomorphisms:
  \begin{align} \label{eq-bigcommutativediagram}
 \xymatrix{
 \cdots  \quad \quad &    \ar[l]_{p_{i_{\ell}}^{i_{\ell-1}}}   \pi_1(M_{i_{\ell}}  , x_{i_{\ell}})    &      &  \ar[ll]_{p_{i_{\ell+1}}^{i_{\ell}}}  \pi_1(M_{i_{\ell+1}} , x_{i_{\ell+1}})       &  \ar[l]_{p_{i_{\ell+2}}^{i_{\ell+1}}}   \quad \cdots            \\
\cdots   \pi_1(M  , y_{j_{\ell-1}})  \ar[d]_{\cong} &   \ar[l]_{\lambda_{\ell}} \pi_1(M_{i_{\ell}}  , x_{i_{\ell}}) \ar[u]_{\cong} &  \ar[l]_{\mu_{\ell}}   \pi_1(M , y_{j_{\ell}})  \ar[d]_{\cong}      &   \ar[l]_{\lambda_{\ell+1}}  \pi_1(M_{i_{\ell+1}} , x_{i_{\ell+1}})  \ar[u]_{\cong}   &  \ar[l]_{\mu_{\ell+1}}  \ar[d]_{\cong}\pi_1(M  , y_{j_{\ell+1}})    \cdots       \\
\cdots    \pi_1(M  , y_{j_{\ell-1}})  &    &  \ar[ll]_{\lambda_{\ell}\circ \mu_{\ell}}   \pi_1(M , y_{j_{\ell}})       &    & \ar[ll]_{\lambda_{\ell+1}\circ \mu_{\ell+1}}   \pi_1(M  , y_{j_{\ell+1}})   \cdots       
} 
\end{align}
Recall that in \eqref{eq-imahes} we defined 
$G_{\ell} = {\rm image}\left\{  p^0_{\ell}  \colon  \pi_1(M_{\ell}, x_{\ell})\longrightarrow G_{0}\right\}$, and the group chain $\cG$ in \eqref{eq-descendingchain} associated to the presentation $\cP$ and basepoint $x_0 \in M_0$.
Define the group chain  
\begin{equation}
\cH = H_{\ell} = {\rm image}\left\{  p^0_{\ell}  \circ \mu_{\ell} \colon  \pi_1(M , y_{j_{\ell}}) \longrightarrow G_{0}\right\} ~, 
\end{equation}
then \eqref{eq-bigcommutativediagram} implies we then obtain a group chain
\begin{equation}\label{eq-interlacedGH}
G_0 \supset G_{i_1} \supset H_{j_1} \supset G_{i_2} \supset \cdots \supset G_{i_{\ell}} \supset H_{j_{\ell}} \supset G_{i_{\ell+1}} \supset H_{j_{\ell +1}} \supset \cdots ~ .
\end{equation}
Set $H = \pi_1(M, y_0)$. As $M$ is connected, for all $\ell > 0$ there is an isomorphism $H_{\ell} \cong H$.
 
Now define  $N_{\ell}$ to be the covering space of $M_0$ associated to the subgroup $H_{\ell}$. The chain of subgroups in \eqref{eq-interlacedGH} then yields a presentation 
$\cP' = \{q_{\ell+1} \colon N_{\ell+1} \to N_{\ell} \mid \ell \geq 1\}$. 
By Theorem~\ref{thm-re},  the spaces $\cS_{\cP}$ and $\cS_{\cP'}$ are homeomorphic. This completes the proof of Theorem~\ref{thm-main1}.

For the proof of Theorem~\ref{thm-main2}, recall that $M$   aspherical means that $\pi_k(M, y_0) = \{0\}$ is the trivial group for all $k > 1$ and choice of basepoint $y_0$.
 The key technical fact  used in the   proof of Theorem~\ref{thm-proM}, the approximate homotopy lifting property in Proposition~\ref{prop-leafonto},  applies equally well for the study of the maps between   higher  homotopy groups in diagram \eqref{eq-commutativediagram3}. Instead of considering lifts of paths to the leaf $L \subset \fM$, one considers lifts of cubes. We require only the conclusion that the following diagram commutes, for each $k > 1$ and $\ell \geq 1$:
 \begin{align} \label{eq-commutativediagrampi1234}
 \xymatrix{
     \pi_k(M_{i_{\ell}}  , x_{i_{\ell}})     &      &  \ar[ll]_{p_{i_{\ell+1}}^{i_{\ell}}}  \pi_k(M_{i_{\ell+1}} , x_{i_1})   \ar[ld]_{\lambda_{\ell+1}}       \\
   &     \pi_k(M , y_{j_{\ell}}) \ar[lu]_{\mu_{\ell}}      &         
} 
\end{align}
The covering map $p_{\ell+1} \colon M_{\ell+1} \to M_{\ell}$ induces an isomorphism on homotopy groups for $k > 1$, hence the horizontal map $p_{i_{\ell+1}}^{i_{\ell}}$ in \eqref{eq-commutativediagrampi1234} is an isomorphism. Thus, $ \pi_k(M_{i_{\ell}}  , x_{i_{\ell}}) = \{0\}$ is trivial for all $\ell \geq 1$, and also $\pi_k(N_{\ell}, n_{\ell}) = \{0\}$ is trivial for the covering space $N_{\ell} \to M_{i_{\ell}}$ defined above, where $n_{\ell} \in N_{\ell}$ is some lift of $x_{i_{\ell}}$. Hence, each $N_{\ell}$ is an aspherical space.  Thus, the lift of $\mu_{\ell}$ to a map $h_{\ell} \colon M \to N_{\ell}$ induces isomorphisms on all homotopy groups, so is a homotopy equivalence. By the assumption that $M$ satisfies the Borel Conjecture, it follows that $h_{\ell}$ is homotopic to a homeomorphism $f_{\ell}$ for each $\ell \geq 1$. Using the homeomorphisms $\{f_{\ell}\}$ we obtain a presentation
$\cP'' = \{f_{\ell}^{-1} \circ q_{\ell+1} \circ f_{\ell+1} \colon M \to M \mid \ell \geq 1\}$
for which  $\cS_{\cP''}$ is homeomorphic to $\cS_{\cP'}$ and hence to $\cS_{\cP}$, which shows the claim of Theorem~\ref{thm-main2}.

\section{Examples}\label{sec-examples}
 
  In this section, we give a   collection of examples, some from the literature and some novel, which illustrate the conclusions of Theorems~\ref{thm-main1} and \ref{thm-main2}.
These examples are based on the basic observation of this work, that the $M$-like hypothesis on a   solenoid is a   shape version of the non co-Hopfian property for manifolds.

 Recall that a group $G$ is \emph{co-Hopfian} if there does not exist an embedding of $G$ to a proper subgroup of itself with finite index, and \emph{non co-Hopfian} otherwise.  (In some literature, a group is said to be co-Hopfian if every embedding of $G$ to a proper subgroup of itself must be an isomorphism, and they refer to our definition as the \emph{finitely co-Hopfian property}. By non  co-Hopfian, we will always mean finitely non co-Hopfian.) The   co-Hopfian concept for groups  was first studied by   Baer in  \cite{Baer1944}, where they are referred to as ``S-groups''. 

A closed connected manifold $M$ is co-Hopfian if every finite covering map $\pi \colon M \to M$ is a homeomorphism, and non co-Hopfian otherwise. Given a proper self-covering map $\pi \colon M \to M$, choose a basepoint $x_1 \in M$ and set $x_0 = \pi(x_1)$. Then $G_0 = \pi_1(M, x_0)$ is a non co-Hopfian group.

Given a non co-Hopfian manifold $M$ with proper self-covering map $\pi \colon M \to M$, we obtain a presentation 
$\cP = \{p_{\ell+1} \colon M_{\ell+1} \to M_{\ell} \mid \ell \geq 0\}$, where each $M_{\ell} = M$ and $p_{\ell+1} = \pi$.  As remarked in Section~\ref{sec-intro}, the associated weak solenoid   $\cS_{\cP}$ is always $M$-like.

Without additional hypotheses on $M$, the property that $G_0 = \pi_1(M, x_0)$ is a non co-Hopfian group does not imply that $M$ is a non co-Hopfian manifold. Here is a simple construction to show that.
 Let $B$ be a non co-Hopfian manifold of dimension $m \geq 4$. 
 Let $M_0$ be the result of attaching the product space  $\mS^2 \times \mS^{m-2}$ to $B$ along the boundary of a  disk in each manifold. Then the fundamental group of $M_0$ equals that of $B$. Note that the second homotopy group $\pi_2(M_0, x_0)$ contains a free abelian summand corresponding to the attached copy of $\mS^2$. For a proper covering $\pi \colon M_1 \to M_0$, the rank of $\pi_2(M_1, x_1)$ is then greater than that of $M_0$, hence they cannot be homeomorphic.  
  
This construction illustrates the requirement that the base manifold $M_0$ be aspherical in Theorem~\ref{thm-main2}. The additional requirement that the coverings of $M_0$ also  satisfy  the Borel Conjecture is made so that all coverings of $M_0$ are then homeomorphic to $M_0$. 
Thus, a   closed aspherical manifold $M$ which satisfies the strongly Borel property \cite[Definition~1.4]{CHL2018a}  is co-Hopfian if and only if its fundamental group is co-Hopfian. 

The non co-Hopfian property was also used in the work of Bi\'s, Hurder and Shive   \cite[Section 3]{BHS2006} to construct many new classes of foliations which generalize the original construction of Hirsch in \cite{Hirsch1975}  that was based on the construction of the Smale solenoid by self-embeddings of a solid torus.

We next recall some   examples of weak solenoids obtained   by explicit constructions of non co-Hopfian closed manifolds, and then discuss some examples obtained from more abstract methods.

 \subsection{Abelian group chains}\label{subsec-abelian}
 
Let $G_0 =\mZ^k$ for $k \geq 1$ be the free abelian group of rank $k$, with basis vectors $\{\vec{e}_1, \ldots , \vec{e}_k\}$. Choose any collection of vectors $\{\vec{v}_1, \ldots, \vec{v}_k\} \subset \mZ^k$ which are linearly independent over $\mQ$. Then define $\phi \colon \mZ^k \to \mZ^k$ by setting $\phi(\vec{e}_i) = \vec{v}_i$ for $1 \leq i \leq k$. The map $\phi$ is an endomorphism of $\mZ^k$, and for most choices we obtain a proper embedding. We obtain a group chain by setting $\cG_{\phi} = \{G_{\ell} \equiv \phi^{\ell}(G_0) \mid \ell \geq 0\}$.

We illustrate this with  the simplest example, where $k=1$. Choose distinct primes $p, q \geq 2$. Define $G_{\ell} = \{n \cdot p^{\ell} \mid n \in \mZ\}$ and $H_{\ell} = \{n \cdot q^{\ell} \mid n \in \mZ\}$, so we obtain  group chains defined by the self-embeddings, 
$\cP = \{G_{\ell} \mid \ell \geq 1\}$ and $\cQ = \{H_{\ell} \mid \ell \geq 1\}$.
Criteria for the equivalence of group chains in $\mZ$ were given by Bing \cite{Bing1960}, McCord \cite{McCord1965}, Aarts and Fokkink \cite{AartsFokkink1991}, as well as   in  \cite[Section 5]{CHL2018a}. 
In particular, $\cP$ and $\cQ$ 
 are not return equivalent,   so we obtain the well-known fact that the solenoids $\cS_{\cP}$ and $\cS_{\cQ}$ are not  homeomorphic.

More generally, every group chain $\cG$  in $\mZ$ defines a scale in the sense of \cite{NekkyPete2011}, see also the discussion in Section \ref{subsec-scale}. 
On the other hand, a monomorphism $\phi \colon \mZ \to \mZ$ is determined by the integer $m = \phi(1)$, which has a finite number of prime divisors, where $m = q_1 \cdots q_k$ with each $q_i$ prime though possibly not distinct. As an example, 
let $\{p_1 < p_2 < \cdots\}$ be an infinite collection of increasing primes. Define $K_{\ell} = \{n \cdot p_1 p_2 \cdots p_{\ell} \mid n \in \mZ\} \cong \mZ$ and let $\cK = \{K_{\ell}  \mid \ell \geq 0\}$ be the corresponding group chain.  Then the group chain  $\cK$ is not return equivalent to any group chain defined by a monomorphism $\phi: \mZ \to \mZ$, and thus the weak solenoid $\cS_{\cK}$ is not homeomorphic to  $\cS(\cG_{\phi})$.

As a consequence, we see that   there are uncountably many circle-like weak solenoids  which are distinct up to homeomorphism, and not homeomorphic to the solenoid obtained from a covering map $\pi \colon \mS^1 \to \mS^1$. 
    Thus, the conclusion of Theorem~\ref{thm-main2} cannot be strengthened in general.

 \subsection{Coverings of the Klein bottle}\label{subsec-RT}
 
  This example is a generalization to arbitrary integer $d>1$ of the example due to Rogers and Tollefson \cite{RT1972a}, see also Fokkink and Oversteegen \cite{FO2002}.  Consider a map of the plane, given by a translation by $\frac{1}{2}$ in the first component, and by reflection in the second component, i.e.
  $$r \times i \colon  \mR^2 \to \mR^2 ~ {\rm where}  ~  (x,y) \mapsto (x+\frac{1}{2},-y).$$
This map commutes with translations by the elements in the integer lattice $\mZ^2 \subset \mR^2$, and so induces the map $r \times i \colon  \mT^2 = \mR^2/\mZ^2 \to \mT^2$ of the torus. This map is an involution, and the quotient space $K = \mT^2/(x,y) \sim r\times i(x,y)$ is homeomorphic to the Klein bottle.

Let  $L \colon  \mT^2 \to \mT^2$ be the  $d$-fold covering map given by  $L(x,y) = (x,dy)$, and form the inverse limit ${\ds \mT_\infty = \lim_{\longleftarrow} \{L \colon  \mT^2 \to \mT^2 \}}$, which  is a solenoid with $2$-dimensional leaves. Let $x_0 = (0,0) \in M_0 = \mT^2$. The fundamental group $G_0 = \mZ^2$ is abelian, so for any $x,y \in \fX_0$ the kernels $K(\cG^x) = K(\cG^y) \cong \mZ$, and every leaf is homeomorphic to an open two-ended cylinder.

The involution $r \times i$ is compatible with the covering maps $L$, and so it induces an involution $(r \times i)_\infty \colon \mT_\infty \to \mT_\infty$, which is seen to have a single fixed point $(0,0, \ldots) \in \mT_\infty$, and permutes the other path-connected components. Let $p \colon  K \to K$ be the $d$-fold covering of  the Klein bottle by itself, given by $p(x,y) = (x,dy)$,  and consider the inverse limit space ${\ds K_\infty = \lim_{\longleftarrow} \{p \colon  K \to K \}}$. Note that taking the quotient by the involution $r \times i$ is compatible with the covering maps $L$ and $p$; that is, $p \circ (r \times i ) =  L$, and so induces the map $i_\infty \colon  \mT_\infty \to K_\infty$ of the inverse limit spaces. Under this map, the path-connected component of the fixed point $(0,0, \ldots)$ is identified so as to become a non-orientable one-ended cylinder. The image of any other path-connected component is an orientable $2$-ended cylinder, so the space $ K_\infty$ is not homogeneous.

Let $x =(x_{\ell})\in K_\infty$ for  $x_{\ell}  \in K$. Then $G_0 = \pi_1(K,x_0) = \langle a,b \mid bab^{-1} = a^{-1} \rangle$. 
The group chain associated to this solenoid is given by the subgroups $G_{\ell} = \langle a^{d^{\ell}}, b   \rangle \subset G_0$. For the group chain $\cG^x$, we have that the kernel $K(\cG^x) = \langle b \rangle $. If a point $y \in K_\infty$ is in a leaf which is orientable, then the kernel of the group chain $\cG^y = \langle b^2 \rangle$. If $y \in K_\infty$ is in the non-orientable leaf, then the kernel $K(\cG^y)$ is conjugate to $K(\cG^x)$.

It follows that  the effective action of $G_0$ on the fiber of the solenoid is that of the dihedral group
$\langle a,b \mid bab^{-1} = a^{-1} ,  b^2 = 1 \rangle$. 
This action is that of the iterated monodromy group associated to a Chebyshev polynomial of degree $d$, as discussed in \cite{Lukina2018b}.

\subsection{Heisenberg group chains}\label{subsec-heisenberg}
 
 Let $\cH$ be the $3$-dimensional Heisenberg Lie group, and let $G_0 = \mH$ be the discrete Heisenberg group. These groups  are presented as the upper-triangular subgroup of the $3 \times 3$   matrices:
 \begin{equation}
\cH \equiv  \left( \begin{array}{ccc} 1 & x &  z\\ 0 & 1 & y\\ 0 & 0 & 1\end{array}\right)  ~ , ~ x,y,z \in \mR \quad , \quad
\mH \equiv  \left( \begin{array}{ccc} 1 & a &  c\\ 0 & 1 & b\\ 0 & 0 & 1\end{array}\right)  ~ , ~ a,b,c \in \mZ
\end{equation}
The group $\cH$ can also be presented  in the form $\{\mR^3, *\}$, with the group operation $*$ given by $(x,y,z)*(x',y',z')=(x+x',y+y',z+z'+xy')$. The group $\mH$  is the   simplest example of a torsion-free nilpotent group which is not abelian. 

Wouter Van Limbeek suggested to the authors the following simple construction of a self-embedding of $G_0$, which yields a group chain that is not normal.
Define  
$\phi \colon G_0 \to G_0$ given by $\phi(a,b,c) = (2a, 2b, 4c)$, which yields a group monomorphism.  
Let 
$$G_{\ell} = \phi^{\ell}(G_0) = \{ (2^{\ell}a, 2^{\ell}b, 4^{\ell} c)\mid a,b,c \in \mZ\}$$
It is then an exercise to check that the normal core of $G_{\ell}$ is the subgroup
$$C_{\ell} =   \{ (4^{\ell}a, 4^{\ell}b, 4^{\ell} c)\mid a,b,c \in \mZ\} \subset G_{\ell} .$$
Note that $G_{2\ell} \subset C_{\ell}$ so the group chain $\{G_{\ell} \mid \ell \geq 0\}$ is equivalent to a normal chain, hence the solenoid defined by self-coverings of $M_0 = \cH/G_0$ will be homogeneous by \cite{FO2002}, and $M$-like by construction.

On the other hand, Jessica Dyer gave in her thesis \cite{Dyer2015}  examples of group chains in $G_0 \subset \mH$ for which the discriminant of the Cantor action constructed from the chain is a Cantor group.  We recall  Example~8.5 from \cite{DHL2016a}.
Choose distinct primes  $p > 1$ and $q>1$, and define
$$H_{\ell} =  \{ (p^{\ell} a, q^{\ell} b , p^{\ell}  c)\mid a,b,c \in \mZ\} ~ .$$
Note that $H_1$ is the image of the set map $\phi(a,b,c) = (pa, qb, pc)$, but this $\phi$ is not a homomorphism for the group law on $\mH$. On the other hand, it is shown in \cite{DHL2016a} that  the solenoid defined by the associated coverings   $M_{\ell} = \cH/H_{\ell}$ of $M_0 = \cH/\mH$ is not homogeneous.

 The following is an interesting   open problem:
 \begin{prob}\label{prob-heisenberg}
 Let $G_0 \subset \mH$ be as above. Does there exists a   subgroup    $G_1 \subset G_0$   of finite index, such that for    $M = \cH/G_1$, there exists an $M$-like weak solenoid which is not homogeneous?
 
  In terms of group chains, 
does there exists $G_1 \subset G_0$  of finite index, and a proper self-embedding $\phi \colon G_1 \to G_1$,  such that the associated group chain $\cG(\phi) \equiv \{G_{\ell} = \phi^{\ell}(G_1) \mid \ell \geq 1\}$ has trivial kernel, and the associated weak solenoid is not homogeneous?
 \end{prob}

 \subsection{Nil-manifolds}\label{subsec-nilmanifolds}

There are many generalizations of the above     constructions  of weak solenoids obtained from lattices in the Heisenberg group. The Heisenberg group $\cH$ is replaced with any connected nilpotent Lie group $\cN$, and $\mH$ is replaced with a cocompact lattice subgroup  $\Gamma \subset \cN$. For example, Belegradek   considered in \cite{Belegradek2003}  when such a lattice must be co-Hopfian, and in particular when they are not. 
Non co-Hopfian subgroups of nilpotent Lie groups were   also studied by Dekimpe, Lee and Potyagailo in \cite{DL2003a,DL2003b,DD2016}, and by Cornulier in \cite{Cornulier2016}.  Here is   a general version of Problem~\ref{prob-heisenberg}:

\begin{prob}\label{prob-nilpotent}
 Let $\cN$ be a connected Lie group. Does there exists a lattice subgroup $\Gamma \subset \cN$ which admits a proper self-embedding $\phi \colon \Gamma \to \Gamma$,  such that the   group chain $\cG(\phi) \equiv \{G_{\ell} = \phi^{\ell}(\Gamma) \mid \ell \geq 0\}$ has trivial kernel, and the associated weak solenoid is not homogeneous.
  \end{prob}

 This problem is closely related to the   works of Reid 
\cite{Reid2014}, Willis \cite{Willis2015}  and Van Limbeek in \cite{vanLimbeek2018,vanLimbeek2017} on the properties of profinite groups derived from a self-embedding of a non co-Hopfian group, in terms of contracting maps for profinite groups.

\subsection{Non co-Hopfian 3-manifolds}\label{subsec-3manifolds}
The question of which compact $3$-manifolds admit proper self-coverings has been   studied in detail by  Gonz{\'a}lez-Acu{\~n}a, Litherland  and  Whitten in the works  \cite{GLiW1994} and \cite{GW1994}, and see also the work \cite{WangWu1994} by Wang and Wu.

\subsection{Higher dimensional constructions}\label{subsec-higher}

Note that the product of any finite collection of  non co-Hopfian manifolds is again non co-Hopfian. Thus, we can form the products of an arbitrary collection of examples as given above to obtain a non co-Hopfian manifold  of arbitrarily large dimension.

 A more general approach is to  consider the problem of which finitely-presented groups $G_0$ are non co-Hopfian? 
One must then also   ask, when such a group can be realized as the fundamental group of a non co-Hopfian manifold? Such constructions usually result in a base manifold of dimension $m \geq 4$. We   recall some results in the literature on when a group is non co-Hopfian.
  
    Endimioni and Robinson give in \cite{ER2005}  some sufficient conditions for a group to be  non co-Hopfian.     
  Delgado and  Timm  consider in \cite{DT2003}   the co-Hopfian condition for the fundamental groups of connected finite complexes, which is certainly a requirement for the construction of non co-Hopfian manifolds.

 Ohshika and Potyagailo gave example in \cite{OP1998}   of a freely indecomposable geometrically finite torsion-free non-elementary Kleinian group which is  non co-Hopfian.   Delzant and  Potyagailo study in \cite{DP2003}      non-elementary geometrically finite Kleinian groups are non co-Hopfian, and Kapovich and Weiss considered the co-Hopf property for word hyperbolic groups in \cite{KW2001}.
 
 On the other hand, the following problem appears to be open:
  \begin{prob}\label{prob-realize}
Let $G_0$ be a finitely-generated non co-Hopfian group. Find conditions on $G_0$ so that it  is realized as the fundamental group of a   non co-Hopfian closed manifold $M_0$.

 \end{prob}

  \subsection{Scale-invariant groups}\label{subsec-scale}
  
   A finitely generated infinite group $G$ is called \emph{scale-invariant} by Nekrashevych  and Pete   \cite{NekkyPete2011},  if there is a group chain $\cG = \{G_{\ell} \mid \ell \geq 0\}$ that each $G_{\ell}$ is    isomorphic to $G$,  and whose kernel  $K(\cG) = \cap_{\ell \geq 0} \ G_{\ell}$ is a finite group.  The assumption that the kernel is finite is a fundamental aspect of this concept. 
  
   Remark~\ref{rmk-kernels} identifies $K(\cG)$ with the fundamental group of the leaf $L_x$ containing the basepoint $x \in \fM$ used to define the group chain. Thus, if $\fM$ is an $M$-like matchbox manifold and $\FfM$ admits a leaf with finite fundamental group,   Theorem~\ref{thm-main1} implies that the fundamental group $G_0 = \pi_1(M, b)$  is scale-invariant. 

A non co-Hopfian group $G_0$ with proper embedding $\phi \colon G_0 \to G_0$  is scale-invariant if the intersection $\cap_{\ell \geq 0} \ \phi^{\ell}(G_0)$ is a finite group.  Nekrashevych and Pete   
addressed in \cite{NekkyPete2011} the problem whether there exists scale invariant groups which are non co-Hopfian. In particular,  they gave examples of scale invariant but non co-Hopfian groups,   obtained from a cross-product of groups, some of which are nilpotent, and some are not.   We briefly recall a special case of this construction.

 Fix $n \geq 1$ and let $G \subset {\rm GL}(n, \mZ)$. Define 
 $$H = \{(\vec{n}, g) \mid  \vec{n} \in \mZ^n ~ , ~ g \in G\}$$
 where the product is defined by 
 $(\vec{n}_1, g_1) \star (\vec{n}_2, g_2) = (\vec{n}_1 + g_1 \cdot \vec{n}_2 , g_1 g_2)$. 
 Then $H$ is nilpotent if and only if the subgroup $G$ is nilpotent. Choose $p > 1$ and define the monomorphism 
 $\phi(\vec{n}_1, g_1) = (p \cdot \vec{n}_1, g_1)$.  The Rogers and Tollefson Example~\ref{subsec-RT} is the simplest example of this construction.
 
 Suppose that $G = \pi_1(B, b)$ for a closed   manifold $B$ which satisfies the strong Borel conjecture. Then let $M$ be the closed manifold which fibers over $B$, with fiber $\mT^n$ twisted by the action of $G$. Then $M$ is again aspherical, and satisfies the strong Borel Conjecture. For example, we can take $G$ to be a  group which acts freely on $\mR^m$ by isometries and let $B = \mR^m/G$ be the quotient. Then $M$ is a non co-Hopfian manifold, and the weak solenoid determined by the map $\phi$ is $M$-like. There are many variations on this construction.


\end{document}